\documentclass[11pt]{article}
\usepackage{mathrsfs}
\usepackage{amssymb}
\usepackage{amsmath}
\usepackage{amsbsy}
\usepackage{epsfig}
\usepackage{enumerate}
\usepackage{bm}
\usepackage{color}

\newcommand{\blue}[1]{\textcolor{blue}{#1}}

\usepackage[colorlinks,linkcolor=blue]{hyperref}

\newcommand{\ARXIV}[1]{\href{http://arXiv.org/abs/#1}{\blue{arXiv:#1}}}

\newtheorem{lemma}{Lemma}[section]
\newtheorem{theorem}{Theorem}[section]
\newtheorem{proposition}{Proposition}[section]

\newtheorem{corollary}{Corollary}[section]
\newtheorem{definition}{Definition}[section]
\newtheorem{remark}{Remark}[section]

\def\pr{\textsf{P}} % the symbol P for probability used the sans serif letter
\def\Sbep{\widehat{\mathbb E}} % the symbol E for sub-linear expectation
\def\cSbep{\widehat{\mathcal E}} % the symbol E for conjugate sub-linear expectation

\def\vSbep{\breve{\mathbb E}}
\def\vcSbep{\breve{\mathcal E}}

\def\Capc{\mathbb V} % the symbol V for capacity under E
\def\cCapc{\mathcal V} % the symbol V for conjugate capacity
\def\upCapc{\widehat{\mathbb V}} % the symbol V with a hat  for a special upper capacity
\def\lowCapc{\widehat{\mathcal V}} % the symbol V with a hat  for a    conjugate capacity of the above
\def\outCapc{\widehat{\mathbb V}^{\ast}}% the symbol  for  the countably additive extension of the upper capacity
\def\outcCapc{\widehat{\mathcal V}^{\ast}} % the symbol  for   for conjugate countably additive capacity

\def\vSbep{\breve{\mathbb E}}

\renewcommand{\baselinestretch}{1.5}

\begin{document}

\begin{center}{\LARGE\bf The sufficient and necessary conditions of the strong law of large numbers under the sub-linear expectations}
\end{center}

\begin{center} {\sc
Li-Xin ZHANG\footnote{This work was Supported by grants from the NSF of China (Grant No.11731012,12031005),   Ten Thousands Talents Plan of Zhejiang Province (Grant No. 2018R52042) and the Fundamental Research Funds for the Central Universities
}
}\\
{\sl \small School  of Mathematical Sciences, Zhejiang University, Hangzhou 310027} \\
(Email:stazlx@zju.edu.cn)\\
%\today
\end{center}

 \renewcommand{\abstractname}{~}
\begin{abstract}
{\bf Abstract:}
In this paper, by establishing a Borel-Cantelli lemma for a capacity which is not necessarily continuous, and a link between a sequence of independent random variables under the sub-linear expectation and a sequence of independent random variables on $\mathbb R^{\infty}$ under a probability,  we give the sufficient and necessary conditions of the strong law of large numbers for independent and identically distributed random variables under the sub-liner expectation, and the  sufficient and necessary conditions for the convergence of an infinite series of independent random variables, without any assumption on the continuity of the capacities.
A purely probabilistic proof of a weak law of large numbers is also given.

{\bf Keywords:}  sub-linear expectation, capacity, strong convergence, law of large numbers

 {\bf AMS 2020 subject classifications:}  60F15, 60F05

% {\bf Short text for running head:}  LLN  under sub-linear expectations
\vspace{-3mm}
\end{abstract}

%%%%%%%%%%%%%%%%%%%%%%%%%%%%%%%%%%%%%%%%%%%%%%%%%%%%%%%%%%%%
\renewcommand{\baselinestretch}{1.2}

%%%%%%%%%%%%%%%%%%%%%%%%%%%%%%%%%%%%%%%%%%%%%%%%%%%%%%%%%%%%
%% Text of article.
%%%%%%%%%%%%%%%%%%%%%%%%%%%%%%%%%%%%%%%%%%%%%%%%%%%%%%%%%%%%
%    Section headings
%\baselineskip 11pt\parindent=10.8pt

\section{ Introduction and notations.}\label{sect1}
\setcounter{equation}{0}

Let $\{X_n;n\ge 1\}$ be a sequence of independent and identically distributed random variables (i.i.d.) on a probability space $(\Omega,\mathcal F, \pr)$. Denote $S_n=\sum_{i=1}^n X_i$. One of the most famous results of probability theory is Kolmogorov (1930)'s strong law of large numbers (c.f., Theorem 3.2.2 of Stout (1974)), which states that
 \begin{equation}\label{eqKol1}
\pr\left(\lim_{n\to \infty}\frac{S_n}{n}=b\right)=1
\end{equation}
if and only if
\begin{equation}\label{eqKol2}
E_{\pr}[|X_1|]<\infty \; \text{ and }\; E_{\pr}[X_1]=b,
\end{equation}
where $E_{\pr}$ is the expectation with respective to the probability measure $\pr$. When the probability measure $\pr$ is uncertain,    one may consider a family $\mathscr{P}$ of probability measures  and define $\Sbep[X]=\sup_{P\in \mathscr{P}}E_P[X]$. Then $\Sbep$ is no longer a linear expectation. It is sub-linear in sense that $\Sbep[aX+bY]\le a\Sbep[X]+b\Sbep[Y]$ if $a,b\ge 0$. Peng (2008,2019) introduced the concepts of independence, identical distribution and $G$-normal random variables under the sub-linear expectation, and established the weak law of large numbers and central limit theorem for independent and identically distributed random variables. Fang et al (2019) obtained the rate of convergence of the weak law of large numbers and central limit theorem.

As for the strong law of large numbers,  Chen (2016) established a Kolmogorov type result. Let $\{X_n; n\ge 1\}$ be a sequence of random variables in a sub-linear expectation  space $(\Omega,\mathscr{H},\Sbep)$ with a related upper capacity $\upCapc$. Chen (2016) showed that, if $\{X_n;n\ge 1\}$ is a sequence of i.i.d. random variables,  the capacity $\Capc$ is continuous, and the following moment condition is satisfied
\begin{equation}\label{eqChen3}
 \Sbep[|X_1|^{1+\alpha}]<\infty \;\;\text{ for some } \alpha>0,
\end{equation}
then
\begin{equation}\label{eqChen1}
\upCapc\left(\liminf_{n\to \infty}\frac{S_n}{n}<-\Sbep[-X_1] \; \text{ and } \;\limsup_{n\to \infty}\frac{S_n}{n}>\Sbep[X_1]\right)=0
\end{equation}
and
\begin{equation}\label{eqChen2}
\upCapc\left(\liminf_{n\to \infty}\frac{S_n}{n}=-\Sbep[-X_1] \right)=1\; \text{ and } \; \upCapc\left(\limsup_{n\to \infty}\frac{S_n}{n}=\Sbep[X_1]\right)=1.
\end{equation}
By establishing the moment inequalities of the maximum partial sums, Zhang (2016) weakened the condition  \eqref{eqChen3}  to
 \begin{equation}\label{eqZhang1}
 C_{\upCapc}(|X_1|):=\int_0^{\infty}\upCapc(|X_1|> x)dx<\infty
 \end{equation}
 and
\begin{equation}\label{eqZhang2}
\Sbep[(|X_1|-c)^+]\to 0  \;\text{ as } c\to \infty.
\end{equation}
The conditions \eqref{eqZhang1}  and \eqref{eqZhang2}  are very close to Kolmogorov's condition \eqref{eqKol2}. Zhang (2016) showed that  \eqref{eqZhang1}  is also a necessary condition. Nevertheless, whether  \eqref{eqZhang2}  is necessary or not is unknown. On the other hand, to make both the direct part and converse part of the Borel-Cantelli lemma are valid for a capacity, it is usually needed to assume that the capacity is continuous  when the strong convergence is considered as in Chen and Hu (2016) and Zhang (2016) etc. However, Zhang (2021b) showed that the assumption of the continuity of a capacity is  very stringent. It is showed that a sub-linear expectation with a continuous capacity is nearly linear.

The purpose of this paper is to obtain the sufficient and necessary conditions for the strong law of large numbers of independent random variables under the sub-linear expectation without the assumption of the continuity of the capacities. In particular it will be shown that, if $\{X_n;n\ge 1\}$ is sequence of i.i.d. random variables in   sub-linear expectation space $(\Omega,\mathscr{H},\Sbep)$ with a regular sub-linear expectation $\Sbep$ and a related upper capacity $\upCapc$ is countably sub-additive (otherwise, $\upCapc$ can be replaced by a countably sub-additive extension), then
$$ \lowCapc\left(\lim_{n\to \infty}\frac{S_n}{n}=b\right)=1 \; \text{ and } b \text{ is finite} $$
if and only if
$$ \text{ (\ref{eqZhang1}) holds and } b=\vSbep[X_1]=\vcSbep[X_1], $$
where $\lowCapc(A)=1-\upCapc(A^c)$, $\vSbep[X]=\lim\limits_{c\to \infty}\Sbep[(-c)\vee X_1\wedge c]$  and $\vcSbep[X]=-\lim\limits_{c\to \infty}\Sbep[(-c)\vee (-X_1)\wedge c]$.

Our main tools are a Borel-Cantelli lemma for a capacity which is not necessarily continuous, and a comparison theorem  for the random variables defined on the product space $\mathbb R^{\infty}$ which gives a link between a sequence of independent random variables on $\mathbb R^{\infty}$   under the sub-linear expectation and a sequence of independent random variables under a probability. By the comparison theorem,  a Kolmogorov's maximal inequality is obtained  and a weak law of large numbers is given with a purely   probabilistic proof.

To state the results, we shall first recall the framework of sub-linear expectation in this section.  We use the framework and notations of Peng (2008, 2019). If one is familiar with these notations, he or she can skip the following paragraphs.  Let  $(\Omega,\mathcal F)$
 be a given measurable space  and let $\mathscr{H}$ be a linear space of real functions
defined on $(\Omega,\mathcal F)$ such that if $X_1,\ldots, X_n \in \mathscr{H}$  then $\varphi(X_1,\ldots,X_n)\in \mathscr{H}$ for each
$\varphi\in C_{l,Lip}(\mathbb R^n)$,  where $C_{l,Lip}(\mathbb R^n)$ denotes the linear space of (local Lipschitz)
functions $\varphi$ satisfying
\begin{eqnarray*} & |\varphi(\bm x) - \varphi(\bm y)| \le  C(1 + |\bm x|^m + |\bm y|^m)|\bm x- \bm y|, \;\; \forall \bm x, \bm y \in \mathbb R^n,&\\
& \text {for some }  C > 0, m \in \mathbb  N \text{ depending on } \varphi. &
\end{eqnarray*}
 We also denote $C_{b,Lip}(\mathbb R^n)$ the space of bounded  Lipschitz
functions.

\begin{definition}\label{def1.1} A  sub-linear expectation $\Sbep$ on $\mathscr{H}$  is a function $\Sbep: \mathscr{H}\to \overline{\mathbb R}$ satisfying the following properties: for all $X, Y \in \mathscr H$, we have
\begin{description}
  \item[\rm (a)]  Monotonicity: If $X \ge  Y$ then $\Sbep [X]\ge \Sbep [Y]$;
\item[\rm (b)] Constant preserving: $\Sbep [c] = c$;
\item[\rm (c)] Sub-additivity: $\Sbep[X+Y]\le \Sbep [X] +\Sbep [Y ]$ whenever $\Sbep [X] +\Sbep [Y ]$ is not of the form $+\infty-\infty$ or $-\infty+\infty$;
\item[\rm (d)] Positive homogeneity: $\Sbep [\lambda X] = \lambda \Sbep  [X]$, $\lambda\ge 0$.
 \end{description}
 Here $\overline{\mathbb R}=[-\infty, \infty]$. The triple $(\Omega, \mathscr{H}, \Sbep)$ is called a sub-linear expectation space. Give a sub-linear expectation $\Sbep $, let us denote the conjugate expectation $\cSbep$of $\Sbep$ by
$$ \cSbep[X]:=-\Sbep[-X], \;\; \forall X\in \mathscr{H}. $$
\end{definition}
By Theorem 1.2.1 of Peng (2019), there exists a family of finite additive linear expectations $E_{\theta}: \mathscr{H}\to \overline{R}$ indexed by $\theta\in \Theta$, such that
\begin{equation}\label{linearexpression} \Sbep[X]=\max_{\theta\in \Theta} E_{\theta}[X] \; \text{ for } X \in \mathscr{H}. \end{equation}
Moreover, for each $X\in \mathscr{H}$, there exists $\theta_X\in \Theta$ such that $\Sbep[X]=E_{\theta_X}[X]$ if $\Sbep[X]$ is finite.

\begin{definition} ({\em See Peng (2008, 2019)})

\begin{description}
  \item[ \rm (i)] ({\em Identical distribution}) Let $\bm X_1$ and $\bm X_2$ be two $n$-dimensional random vectors defined
respectively in sub-linear expectation spaces $(\Omega_1, \mathscr{H}_1, \Sbep_1)$
  and $(\Omega_2, \mathscr{H}_2, \Sbep_2)$. They are called identically distributed, denoted by $\bm X_1\overset{d}= \bm X_2$  if
$$ \Sbep_1[\varphi(\bm X_1)]=\Sbep_2[\varphi(\bm X_2)], \;\; \forall \varphi\in C_{b,Lip}(\mathbb R^n). $$
 A sequence $\{X_n;n\ge 1\}$ of random variables is said to be identically distributed if $X_i\overset{d}= X_1$ for each $i\ge 1$.
\item[\rm (ii)] ({\em Independence})   In a sub-linear expectation space  $(\Omega, \mathscr{H}, \Sbep)$, a random vector $\bm Y =
(Y_1, \ldots, Y_n)$, $Y_i \in \mathscr{H}$ is said to be independent to another random vector $\bm X =
(X_1, \ldots, X_m)$ , $X_i \in \mathscr{H}$ under $\Sbep$  if for each test function $\varphi\in C_{l,Lip}(\mathbb R^m \times \mathbb R^n)$
we have
$ \Sbep [\varphi(\bm X, \bm Y )] = \Sbep \big[\Sbep[\varphi(\bm x, \bm Y )]\big|_{\bm x=\bm X}\big],$
whenever $\overline{\varphi}(\bm x):=\Sbep\left[|\varphi(\bm x, \bm Y )|\right]<\infty$ for all $\bm x$ and
 $\Sbep\left[|\overline{\varphi}(\bm X)|\right]<\infty$.

 A sequence of random variables $\{X_n; n\ge 1\}$
 is said to be independent, if  $X_{i+1}$ is independent to $(X_{1},\ldots, X_i)$ for each $i\ge 1$.
 \end{description}
\end{definition}

Next, we consider the capacities corresponding to the sub-linear expectations. Let $\mathcal G\subset\mathcal F$. A function $V:\mathcal G\to [0,1]$ is called a capacity if
$$ V(\emptyset)=0, \;V(\Omega)=1 \; \text{ and } V(A)\le V(B)\;\; \forall\; A\subset B, \; A,B\in \mathcal G. $$
It is called to be sub-additive if $V(A\bigcup B)\le V(A)+V(B)$ for all $A,B\in \mathcal G$  with $A\bigcup B\in \mathcal G$.

 Let $(\Omega, \mathscr{H}, \Sbep)$ be a sub-linear expectation space.  We denote   $(\upCapc,\lowCapc)$ be a pair of  capacities by
 \begin{equation}\label{eq1.3} \upCapc(A):=\inf\{\Sbep[\xi]: I_A\le \xi, \xi\in\mathscr{H}\},\; \lowCapc(A)=1-\upCapc(A^c),  \forall A\in \mathcal F,
\end{equation}
where $A^c$ is the complement set of $A$. Then $\upCapc$ is a sub-additive capacity with the property that
\begin{equation}\label{eq1.4}
   \Sbep[f]\le \upCapc(A)\le \Sbep[g]\;\;
\text{ if } 0\le f\le I_A\le g, f,g \in \mathscr{H} \text{ and } A\in \mathcal F.
\end{equation}
We call $\upCapc$ and $\lowCapc$ the upper and the lower capacity, respectively.

 Also, we define the  Choquet integrals/expecations $(C_{\upCapc},C_{\lowCapc})$  by
$$ C_V[X]=\int_0^{\infty} V(X\ge t)dt +\int_{-\infty}^0\left[V(X\ge t)-1\right]dt $$
with $V$ being replaced by $\upCapc$ and $\lowCapc$ respectively.
If $\Capc$ on the sub-linear expectation space $(\Omega,\mathscr{H},\Sbep)$ and $\widetilde{\Capc}$ on the sub-linear expectation space $(\widetilde{\Omega},\widetilde{\mathscr{H}},\widetilde{\mathbb E})$  are two capacities having the property \eqref{eq1.4}, then for any random variables $X\in \mathscr{H}$ and $\tilde X\in \widetilde{\mathscr{H}}$ with $X\overset{d}=\tilde X$, we have
\begin{equation}\label{eqV-V}
\Capc(X\ge x+\epsilon)\le \widetilde{\Capc}(\tilde X\ge x)\le \Capc(X\ge x-\epsilon)\;\; \text{ for all } \epsilon>0 \text{ and } x.
\end{equation}
In fact, let $f\in C_{b,Lip}(\mathbb R)$ such that $I\{y\ge x+\epsilon\}\le f(y)\le I\{y\ge x\}$.
Then
$$\Capc(X\ge x+\epsilon)\le \Sbep[f(X)]=\widetilde{\mathbb E}[f(\tilde X)]\le \widetilde{\Capc}(X\ge x), $$
and similar $\widetilde{\Capc}(\tilde X\ge x+\epsilon)\le \Capc(X\ge x)$.  From \eqref{eqV-V}, it follows that $\Capc(X\ge x)=\widetilde{\Capc}(\tilde X\ge x)$ if $x$ is a continuous point of both functions $\Capc(X\ge y)$ and $\widetilde{\Capc}(\tilde X\ge y)$. Since, a monotone function has at most countable number of discontinuous points. So
$$
\Capc(X\ge x)=\widetilde{\Capc}(\tilde X\ge x)\;\;  \text{ for all but except countable many }  x,
$$
and then
\begin{equation}\label{eqV-V3} C_{\Capc}(X)=C_{\widetilde{\Capc}}(\tilde X).
\end{equation}

Because a capacity $\upCapc$ may be not countably sub-additive so that the Borel-Cantelli lemma is not valid,  we  consider its countably sub-additive extension   $\outCapc$ which defined by
\begin{equation}\label{outcapc} \outCapc(A)=\inf\Big\{\sum_{n=1}^{\infty}\upCapc(A_n): A\subset \bigcup_{n=1}^{\infty}A_n\Big\},\;\; \outcCapc(A)=1-\outCapc(A^c),\;\;\; A\in\mathcal F.
\end{equation}
 As shown in Zhang (2016), $\outCapc$ is countably sub-additive, and $\outCapc(A)\le \upCapc(A)$. Further, $\upCapc$ (resp. $\outCapc$) is the largest sub-additive (resp. countably sub-additive) capacity in sense that  if $V$ is also a sub-additive (resp. countably sub-additive) capacity satisfying $V(A)\le \Sbep[g]$ whenever $I_A\le g\in \mathscr{H}$, then $V(A)\le \Capc(A)$ (resp. $V(A)\le \outCapc(A)$).

 Besides $\outCapc$, another countably sub-additive capacity generated by $\Sbep$ can be defined as follows:
 \begin{equation}\label{tildecapc} \mathbb C^{\ast}(A)=\inf\Big\{\lim_{n\to\infty}\Sbep[\sum_{i=1}^n g_i]: I_A\le \sum_{n=1}^{\infty}g_n, 0\le g_n\in\mathscr{H}\Big\},\;\;\; A\in\mathcal F.
\end{equation}
It can be shown that the out capacity $c^{\prime}$ defined in Example 6.5.1 of Peng (2019) coincides with  $\mathbb C^{\ast}$ if $\mathscr{H}$ is chosen as the family of (bounded) continuous functions on a metric space $\Omega$.

For real numbers $x$ and $y$, denote $x\vee y=\max(x,y)$, $x\wedge y=\min(x,y)$, $x^+=\max(0,x)$ and $x^-=\max(0,-x)$. For a random variable $X$, because $XI\{|X|\le c\}$  may be not in $\mathscr{H}$, we will truncate it in the form $(-c)\vee X\wedge c$ denoted by $X^{(c)}$, and define
$\vSbep[X]=\lim\limits_{c\to\infty}\Sbep[X^{(c)}]$ if the limit exists, and $\vcSbep[X]=-\vSbep[-X]$.
\begin{proposition}\label{prop1.1} Consider a subspace of $\mathscr{H}$ as
\begin{equation}\mathscr{H}_1=\big\{X\in\mathscr{H}: \lim_{c,d\to \infty}\Sbep\big[(|X|\wedge d-c)^+\big]=0\big\}.
\end{equation}
Then for any $X\in\mathscr{H}_1$, $\vSbep[X]$ is well defined, and $(\Omega,\mathscr{H}_1,\vSbep)$ is a sub-linear space.
\end{proposition}
{\bf Proof}.    For any $X\in\mathscr{H}_1$ and $0<c_1,c_2\le d$  we have
\begin{align*}  \Sbep[|(-c_1)\vee X  \wedge c_2 -X^{(d)}|]
\le   \Sbep\left[\big((|X|\wedge d-(c_1\wedge c_2)\big)^+\right].
\end{align*}
Hence
$$ \left|\Sbep[ X^{(c)}]-\Sbep[X^{(d)}]\right|\to 0 \;\; \text {as } c,d\to \infty, $$
which implies that $\vSbep[X]=\lim\limits_{c\to \infty}\Sbep[ X^{(c)}]$ exists and is finite. Further,
\begin{equation}\label{eqproposition1.1}\lim_{c_1,c_2\to \infty}  \Sbep[|(-c_1)\vee X \wedge c_2 ]=\vSbep[X].
\end{equation}
Note $(\lambda X)^{(c)}=X^{(c/\lambda)}$ for $\lambda>0$. It is obvious  that $\vSbep[X]=\lambda\vSbep[X]$ for $\lambda>0$.
Finally, for any $X,Y\in \mathscr{H}_1$ and $c>0$, we have $X+Y\in\mathscr{H}_1$ and
\begin{align*}   \big(X+Y\big)^{(c)}
\le   (-c/2)\vee X \wedge (3c/2)+(-c/2)\vee Y \wedge (3c/2).
\end{align*}
By \eqref{eqproposition1.1}, $\vSbep[X+Y]\le \vSbep[X]+\vSbep[Y]$. The monotonicity and constant preserving  for $\vSbep$ are obvious.  The proof is completed.  $\Box$

Let
\begin{equation}\label{linearexpression2.1}\mathscr{E}=\{E: \mathscr{H}_1\to \mathbb R \text{ is a finite additive linear expectation with } E\le \vSbep\}.
\end{equation}
By Theorem 1.2.1 of Peng (2019),
\begin{equation}\label{linearexpression2.2} \vSbep[X]=\max_{E\in\mathscr{E}} E[X] \; \text{ for } X \in \mathscr{H}_1, \end{equation}
and moreover, for each $X\in \mathscr{H}_1$, there exists $E\in \mathscr{E}$ such that $\vSbep[X]=E[X]$. For the vector $\bm X=(X_1,\ldots,X_d)$, we denote
$ \vSbep[\bm X]=(\vSbep[X_1],\ldots,\vSbep[X_d])$, $ \Sbep[\bm X]=(\Sbep[X_1],\ldots,\Sbep[X_d])$ and $E[\bm X] =(E[X_1],\ldots,E[X_d])$ for $E\in \mathscr{E}$.

Finally, a random variable $X$ is called tight (under a capacity $\Capc$ satisfying \eqref{eq1.4} if $\Capc(|X|\ge c)\to 0$ as $c\to\infty$. It is obvious that if $\Sbep[|X|]<\infty$, or $\vSbep[|X|]<\infty$ or $C_{\Capc}(|X|)<\infty$, then $X$ is tight.

%%%%%%%%%%%%%%%%%%%%%%%%%%%%%%%%%%%%%%%%%%%%%%%%%%%

  \section{Basic Tools}\label{sectTools}
\setcounter{equation}{0}

In this section, we give some results which are basic tools for   establishing  the law of large numbers as well as other limit theorems. The first one gives a link between the capacity and a probability measure.

\begin{proposition}\label{lem3.4}  Let $(\Omega,\mathscr{H},\Sbep)$ be a sub-linear expectation space with a capacity $\Capc$ satisfying \eqref{eq1.4},  and $\{X_n; n\ge 1\}$  be a sequence of  random variables in   $(\Omega,\mathscr{H},\Sbep)$. We can find a new sub-linear space $(\widetilde{\Omega},\widetilde{\mathscr{H}},\widetilde{\mathbb E})$ defined on a metric space $\widetilde{\Omega}=\mathbb R^{\infty}$, with  a sequence $\{\tilde X_1,\tilde X_2,\ldots\}$  of  random variables  and a set function $\widetilde{V}:\widetilde{\mathcal F}\to [0,1]$ on it satisfying the following properties, where $\widetilde{\mathcal F}=\sigma(\widetilde{\mathscr{H}})$.
\begin{description}
  \item[\rm (a) ]  $(X_1,X_2,\ldots,X_n)\overset{d}=(\tilde X_1,\tilde X_2,\ldots,\tilde X_n)$, $n=1,2,\ldots$, i.e.,
  $$\widetilde{\mathbb E}[\varphi(\tilde{X}_1,\ldots,\tilde{X}_n)]=  \Sbep[\varphi(X_1,\ldots,X_n)], \;\;   \varphi\in C_{l,Lip}(\mathbb R^n), n\ge 1. $$
  In particular, if $\{X_n;\ge 1\}$ are independent under $\Sbep$, then $\{\tilde X_n;n\ge 1\}$ are independent under $\widetilde{\mathbb E}$.
  \item[\rm (b)] Define
$$ \widetilde{V}(A)=\sup_{P\in\widetilde{\mathscr{P}}}P(A),\;\;  A\in \widetilde{\mathcal F},$$
where $\widetilde{\mathscr{P}}$ is the family of all probability measures $P$ on $(\widetilde{\Omega},\widetilde{\mathcal F})$ with the property
$$ P[\varphi]\le \widetilde{\mathbb E}[\varphi] \; \text{ for bounded }\varphi\in  \widetilde{\mathscr{H}},  $$
and $\widetilde{V}\equiv 0$ if $\widetilde{\mathscr{P}}$ is empty.
Then $\widetilde{V}:  \widetilde{\mathcal F}\to [0,1]$ is a countably sub-additive and nondecreasing function, and   $\widetilde{V}\le  \widetilde{\mathbb C}^{\ast}\le \widetilde{\Capc}^{\ast}\le \widetilde{\Capc}$,
  where $\widetilde{\Capc}$,  $\widetilde{\Capc}^{\ast}$ and $\widetilde{\mathbb C}^{\ast}$ are defined on $(\widetilde{\Omega},\widetilde{\mathscr{H}},\widetilde{\mathbb E})$ in the same way as $\upCapc$,  $\outCapc$  and $\mathbb C^{\ast}$ on $(\Omega,\mathscr{H},\Sbep)$, respectively.

Here and in the sequel, for a probability measure $P$ and a measurable function $X$, $PX$ is defined to be the expectation $\int X dP$.
  \item[\rm (c) ]  If   each $X_n$ is tight, then $\widetilde{\mathscr{P}}$ is a weakly compact family of probability measures on $\widetilde{\Omega}$,
\begin{equation}\label{eqexpression} \widetilde{\mathbb E}[\varphi]=\sup_{P\in\widetilde{\mathscr{P}}}P[\varphi]  \; \text{ for bounded }\varphi\in  \widetilde{\mathscr{H}},
\end{equation}
and $\widetilde{V}$ is a countably sub-additive capacity,
  \begin{equation}\label{eq1.4V}
   \widetilde{\mathbb E}[f]\le \widetilde{V}(A)\le \widetilde{\mathbb E}[g]\;\;
\text{ if } 0\le f\le I_A\le g, f,g \in \widetilde{\mathscr{H}} \text{ and } A\in \widetilde{\mathcal F}.
\end{equation}
  \item[\rm (d) ] If $\{X_n;\ge 1\}$ are independent under $\Sbep$ and each $X_n$ is tight, then for any sequence of vectors $\{\bm \xi_k=(X_{n_{k-1}+1},\ldots,X_{n_k}); k\ge 1\}$  and a sequence  $\{E_k;k\ge 1\}$   of finite additive linear expectations on $\mathscr{H}_b=\{f\in \mathscr{H}; f \text{ is bounded}\}$ with $E_k\le \Sbep$, where $1=n_0<n_1<n_2<\ldots$, there exists a probability measure $Q$ on $\sigma(\tilde{X}_1,\tilde{X}_2,\ldots)$ such that $\{\tilde{\bm \xi}_k=(\tilde X_{n_{k-1}+1},\ldots,\tilde X_{n_k});k\ge 1\}$ is a sequence of independent random vectors under $Q$,
      \begin{equation}\label{eqP<V1}
Q\left[\varphi(\tilde{\bm \xi}_k)\right]=E_k\left[\varphi(\bm \xi_k)\right]\; \text{ for all } \varphi\in C_{b,Lip}(\mathbb R^{n_k-n_{k-1}}),
\end{equation}
\begin{equation}\label{eqP<V2}
Q\left[\varphi(\tilde{X}_1,\ldots,\tilde{X}_n)\right]\le \Sbep\left[\varphi(X_1,\ldots,X_d)\right]\; \text{ for all } \varphi\in C_{b,Lip}(\mathbb R^{n})
\end{equation}
and
\begin{align}\label{eqP<V3}
\widetilde{v}\left((\tilde{X}_1,\tilde{X}_2,\ldots)\in B\right)&\le Q\left((\tilde{X}_1,\tilde{X}_2,\ldots)\in B\right)\le \widetilde{V}\left((\tilde{X}_1,\tilde{X}_2,\ldots)\in B\right) \\
& \; \text{ for all } B\in \mathscr{B}(\mathbb R^{\infty}),\nonumber
\end{align}
where $\widetilde{v}=1-\widetilde{V}$.
\end{description}
\end{proposition}

\begin{remark} When $X_1,X_2,\ldots$ are bounded random variables, then \eqref{eqP<V1} and \eqref{eqP<V2} hold for all $\varphi\in C_{l,Lip}$. When $\bm X_1,\bm X_2,\ldots$ are multi-dimensional random vectors, Proposition \ref{lem3.4} remains true.
\end{remark}

{\bf Proof}. This lemma is proved in Zhang (2021b). We summarize the results and  the proof here for the convenience of reading  and the completeness of this paper.  We  use the the key idea in   Lemma 1.3.5 of Peng (2019) to construct the new sub-linear expectation in the real space.  Let $\widetilde{\Omega}=\mathbb R^{\infty}$, $\widetilde{\mathcal{F}}=\mathscr{B}(\mathbb R^{\infty})$ and
$$\widetilde{\mathscr{H}}=\left\{\varphi(x_1,\ldots,x_n): \varphi\in C_{l,Lip}(\mathbb R^d), d\ge 1, \text{ for } \bm x=(x_1,x_2,\ldots)\in \widetilde{\Omega}\right\}. $$
Define
$$ \widetilde{\mathbb E}[\varphi]=\Sbep[\varphi(X_1,\ldots,X_n)], \;\; \varphi \in C_{l,Lip}(\mathbb R^d). $$
Then $\widetilde{\mathbb E}$ is a sub-linear expectation on $(\widetilde{\Omega},\widetilde{\mathscr{H}})$. Define the random variable $\tilde{X}_i$ by
$\tilde{X}_i(\tilde{\omega})=x_i$ for $\tilde{\omega}=(x_1,x_2,\ldots)\in \widetilde{\Omega}$. Then
$$\widetilde{\mathbb E}[\varphi(\tilde{X}_1,\ldots,\tilde{X}_n)]= \widetilde{\mathbb E}[\varphi]=\Sbep[\varphi(X_1,\ldots,X_n)], \;\;   \varphi\in C_{l,Lip}(\mathbb R^d). $$
It follows that $(\tilde{X}_1,\ldots,\tilde{X}_n)\overset{d}=(X_1,\ldots,X_n)$ for $n=1,2\ldots$. (a) is proved,  and  (b) is obvious.

For (c),  suppose that each $X_n$ is tight.
For the new sub-linear expectation, we also have the expression \eqref{linearexpression}:
$$ \widetilde{\mathbb E}[\tilde{X}]=\max_{\widetilde{\theta}\in \widetilde{\Theta}} E_{\theta}[\tilde{X}] \; \text{ for } \tilde{X} \in \widetilde{\mathscr{H}},$$
 for  a family of finite additive linear expectations $E_{\widetilde{\theta}}: \widetilde{\mathscr{H}}\to \overline{\mathbb R}$ indexed by $\widetilde{\theta}\in \widetilde{\Theta}$, and for each $\tilde{X}\in \widetilde{\mathscr{H}}$, there exists $\widetilde{\theta}_{\tilde{X}}\in \widetilde{\Theta}$ such that $\widetilde{\mathbb E}[\tilde{X}]=E_{\theta_{\tilde{X}}}[\tilde X]$ if $\widetilde{\mathbb E}[\tilde{X}]$ is finite. For each $E_{\widetilde{\theta}}$, consider the finite additive linear expectation $E_{\widetilde{\theta}}$ on $C_{b,Lip}(\mathbb R^p)$. For any sequence $C_{b,Lip}(\mathbb R^p)\ni \phi_n\searrow 0$, we have $\sup_{|\bm x|\le c}|\varphi_n(\bm x)|\to 0$, and so
  $$ E_{\widetilde{\theta}}[\varphi_n]\le \Sbep \left[\varphi_n(X_1,\ldots,X_d)\right]\le \sup_{|\bm x|\le c}|\varphi_n(\bm x)| +\sum_{j=1}^p\|\varphi_1\|\Capc(|X_j|>c)\to 0$$
  as $n\to \infty$ and then $c\to \infty$, by the tightness of $X_j$, where $\|\varphi\|=\sup_{\bm x}|\varphi(\bm x)|$.  Then, as shown in Lemma 1.3.5 of Peng (2019), by Daniell-Stone's theorem, there exists a family of probability measures $P_{\widetilde{\theta},p}$  on $(\mathbb R^p,\mathscr{B}(\mathbb R^p))$ such that
 $$ E_{\widetilde{\theta}}[\varphi]=P_{\widetilde{\theta},p}[\varphi]=\int \varphi(x_1,\ldots,x_p) P_{\widetilde{\theta},p}(d x_1,\ldots,d x_p), \; \varphi\in C_{b,Lip}(\mathbb R^p).  $$
It is obvious that $\{P_{\widetilde{\theta},p};p\ge 1\}$ is a Kolmogorov's consistency system. By Kolmogorov's existence theorem, there is a unique probability measure $P$  on $(\mathbb R^{\infty},\mathscr{B}(\mathbb R^{\infty}))$ such that $P_{\theta}\big|_{\mathscr{B}(\mathbb R^p)}=P_{\widetilde{\theta},p}$. Hence
$$P_{\widetilde{\theta}}[\varphi]= E_{\widetilde{\theta}}[\varphi]\le \widetilde{\mathbb E}[\varphi], \; \varphi\in C_{b,Lip}(\mathbb R^p).  $$
Recall that $\widetilde{\mathscr{P}}$ is the family of all probability measures $P$ on $(\mathbb R^{\infty},\mathscr{B}(\mathbb R^{\infty}))$ with the property
$$ P[\varphi]\le \widetilde{\mathbb E}[\varphi], \; \text{ for bounded }\varphi\in  \widetilde{\mathscr{H}}.  $$
Then for any bounded $\varphi\in  \widetilde{\mathscr{H}}$,
$$ \widetilde{\mathbb E}[\varphi]=\sup_{\tilde{\theta}\in \widetilde{\Theta}}E_{\tilde{\theta}}[\varphi]=\sup_{\tilde{\theta}\in \widetilde{\Theta}}P_{\tilde{\theta}}[\varphi]\le \sup_{P\in\widetilde{\mathscr{P}}}P[\varphi]\le  \widetilde{\mathbb E}[\varphi].$$
It follows that \eqref{eqexpression} holds
and for each bounded $\varphi\in  \widetilde{\mathscr{H}}$ there exists a $P\in\widetilde{\mathscr{P}}$ such that $P[\varphi]= \widetilde{\mathbb E}[\varphi]$.

Suppose $0\le f\le I_A\le g$,  $f(\bm x)=f(x_1,\ldots,x_p)$, $g(\bm x)=g(x_1,\ldots,x_p) \in \widetilde{\mathscr{H}}$   and $A\in \widetilde{\mathcal F}$. Then
$$ P[f]\le P(A)\le P[g\wedge 1]. $$
By \eqref{eqexpression} and taking the supremum over $P\in\widetilde{\mathscr{P}}$, it follows that
$$ \Sbep[f(X_1,\ldots,X_p)]=\widetilde{\mathbb E}[f]\le \widetilde{V}(A)\le \widetilde{\mathbb E}[g\wedge 1]\le \widetilde{\mathbb E}[g]=\Sbep[g(X_1,\ldots,X_p)]. $$
\eqref{eq1.4V} is proved.
At last, we show that $\widetilde{\mathscr{P}}$ is weakly compact. For any $\epsilon>0$, by the tightness of $X_i$, there exists a constant $C_i$ such  that
$\Capc(|X_i|\ge C_i)<\epsilon/2^i$. Then  $\widetilde{V}(\bm x:|x_i|\ge 2C_i)\le \Capc(|X_i|\ge C_i)<\epsilon/2^i$. Let $K=\bigotimes_{i=1}^{\infty}[-2C_i,2C_i]$. Then $K$ is a compact subset in the space $\mathbb R^{\infty}$ with a metric defined by $d(\bm x,\bm y)=\sum_{i=1}^{\infty}(|x_i-y_i|\wedge 1)/2^i$. Note
$$ \widetilde{V}(\bm x\not\in K)\le \sum_{i=1}^{\infty}\widetilde{V}(\bm x:|x_i|\ge 2C_i)\le \sum_{i=1}^{\infty} \epsilon/2^i<\epsilon. $$
It follows that $\widetilde{\mathscr{P}}$ is tight and so is relatively weakly compact. Assume $\widetilde{\mathscr{P}}\ni P_n\implies P$. It is obvious that
$$ P[f]=\lim_{n\to \infty}P_n[f]\le \Sbep[f] \; \text{ for bounded } f\in \widetilde{\mathscr{H}}. $$
Hence $P\in \widetilde{\mathscr{P}}$. It follows that $\widetilde{\mathscr{P}}$ is close and so is weakly compact.
 (c) is proved.

Now, we show (d).   Consider the  linear operator $\widetilde{E}_k$ on $C_{b,Lip}(\mathbb R^{n_k-n_{k-1}})$ defined by
$$ \widetilde{E}_k[\varphi]=E_k \left[\varphi(\bm \xi_k)\right], \;\;\varphi\in C_{b,Lip}(\mathbb R^{n_k-n_{k-1}}). $$
Then
$$\widetilde{E}_k[\varphi] \le \Sbep \left[\varphi(\bm \xi_k)\right], \;\;\varphi\in C_{b,Lip}(\mathbb R^{n_k-n_{k-1}}). $$
  If $C_{l,Lip}(\mathbb R^{n_k-n_{k-1}})\ni \varphi_n\searrow 0$, then $\sup_{|\bm x|\le c}|\varphi_n(\bm x)|\to 0$ and
  $$ \widetilde{E}_k[\varphi_n]\le \Sbep \left[\varphi_n(\bm \xi_k)\right]\le \sup_{|\bm x|\le c}|\varphi_n(\bm x)| +\|\varphi_1\|\Capc(|\bm \xi_k|>c)\to 0$$
  as $n\to \infty$ and then $c\to \infty$, where $\|\varphi\|=\sup_{\bm x}|\varphi(\bm x)|$.
By Daniell-Stone's theorem again,   there exists a probability $Q_k$ on $\mathbb R^{n_k-n_{k-1}}$ such that
$$  Q_k[\varphi]=\widetilde{E}_k[\varphi]\le \Sbep[\varphi(\bm \xi_k)],\;\;\forall \varphi\in C_{b,Lip}(\mathbb R^{n_k-n_{k-1}}).
$$
Now, we introduce  a product probability measure on $\mathbb R^{\infty}$ defined by
$$ Q=Q_1\big|_{\mathbb R^{n_1}}\times Q_2\big|_{\mathbb R^{n_2-n_1}}\times \cdots. $$
Then, under the probability measure $Q$,  for any $A_i\in \mathscr{B}(\mathbb R^{n_i-n_{i-1}})$, $i=1,\ldots,d$, $d\ge 1$,
$$ Q(A_1\times  \cdots \times A_d)=Q(A_1)\cdots Q(A_d), $$
That is
$$ Q(\tilde{\bm \xi}_1\in A_1,  \cdots,\tilde{\bm \xi}_d\in A_d)=Q(\tilde{\bm \xi}_1\in A_1)\cdots Q(\tilde{\bm \xi}_d\in A_d). $$
So, $\tilde{\bm \xi}_1,\tilde{\bm \xi}_2,\cdots$ is a sequence of independent random variables under $Q$. Further,
\begin{equation}\label{eqprooflem3.4.1}
Q[\varphi(\tilde{\bm \xi}_k)]=Q_k[\varphi)]=\widetilde{E}_k[\varphi]=E_k[\varphi(\bm \xi_k)]\le \Sbep[\varphi(\bm \xi_k)],\;\;\forall \varphi\in C_{b,Lip}(\mathbb R^{n_k-n_{k-1}}).
\end{equation}
 \eqref{eqP<V1} is proved.

Note that  for every $\varphi(\bm z_1,\ldots,\bm z_d)\in C_{b,Lip}(\mathbb R^{n_d})$, where $\bm z_i=(x_{n_{i-1}+1},\ldots,x_{n_i})$,
$$ Q[\varphi(\bm z_1,\ldots,\bm z_{d-1},\tilde{\bm \xi}_d)]=Q_d[\varphi(\bm z_1,\ldots,\bm z_{d-1},\cdot)]
 \le \Sbep[\varphi(\bm z_1,\ldots,\bm z_{d-1},\tilde{\bm \xi}_d)]  $$
by \eqref{eqprooflem3.4.1}. Write the functions of $(\bm z_1,\ldots,\bm z_{d-1})$  in the left hand and right hand by $\varphi_1(\bm z_1,\ldots,\bm z_{d-1})$   and
 $\varphi_2(\bm z_1,\ldots,\bm z_{d-1})$, respectively. Note that $\tilde{\bm \xi}_1, \ldots,   \tilde{\bm \xi}_d$ are independent under both $Q$ and $\widetilde{\mathbb E}$, and  $\bm \xi_1,\ldots,\bm \xi_d$ are independent under $\Sbep$. We have that
\begin{align*}
 &Q[\varphi(\bm z_1,\ldots,\bm z_{d-2}, \tilde{\bm \xi}_{d-1}, \tilde{\bm \xi}_d)]
 =Q\left[\varphi_1(\bm z_1,\ldots,\bm z_{d-2},\tilde{\bm \xi}_{d-1})\right]\\
 \le &Q\left[\varphi_2(\bm z_1,\ldots,\bm z_{d-2},\tilde{\bm \xi}_{d-1})\right]
 \le \Sbep\left[\varphi_2(\bm z_1,\ldots,\bm z_{d-2}, \bm \xi_{d-1})\right] \\
 =&\Sbep[\varphi(\bm z_1,\ldots,\bm z_{d-2}, \bm \xi_{d-1},  \bm \xi_d)]=\widetilde{\mathbb E}[\varphi(\bm z_1,\ldots,\bm z_{d-2}, \tilde{\bm \xi}_{d-1}, \tilde{\bm \xi}_d)],
 \end{align*}
by \eqref{eqprooflem3.4.1} again. By induction, we conclude  that
$$ Q[\varphi(\tilde{\bm \xi}_1, \ldots,   \tilde{\bm \xi}_d)]\le \widetilde{\mathbb E}[\varphi(\tilde{\bm \xi}_1,\ldots,\tilde{\bm \xi}_d)], \;\text{ for all } \varphi\in C_{b,Lip}(\mathbb R^{n_d}), \; d\ge 1. $$
Now, for each $\varphi\in C_{b,Lip}(\mathbb R^{n})$, $\varphi\circ\pi_{n_k\to n}$ is also a function in $ C_{b,Lip}(\mathbb R^{n_d})$ when $n\le n_d$, where $\pi_{n_k\to n}:\mathbb R^{n_k}\to \mathbb R^n$ is the projection map. It follows that
\begin{align*}
Q[\varphi]=& Q[\varphi(\tilde{X}_1,\ldots,\tilde{X}_n)]=Q[\varphi\circ\pi_{n_k\to n}(\tilde{\bm \xi}_1, \ldots, \tilde{\bm \xi}_d)]\\
\le & \widetilde{\mathbb E}[\varphi\circ\pi_{n_k\to n}(\tilde{\bm \xi}_1, \ldots, \tilde{\bm \xi}_d]=\widetilde{\mathbb E}[\varphi(\tilde{X}_1,\ldots,\tilde{X}_n)]  \\
= & \Sbep[\varphi(X_1,\ldots,X_n)]=\widetilde{\mathbb E}[\varphi].
\end{align*}
That is, $Q[\varphi]\le \widetilde{\mathbb E}[\varphi]$ for all bounded $\varphi\in \mathscr{H}$.
Hence,   $Q\in \widetilde{\mathscr{P}}$ and \eqref{eqP<V2} holds. So, for each $B\in\mathscr{B}(\mathbb R^{\infty})$,
$$  Q\left((\tilde X_1, \tilde X_2,\ldots)\in B\right)=Q(B)\le \widetilde{V}(B)=\widetilde{V}\left((\tilde X_1, \tilde X_2,\ldots)\in B\right), $$
by the definition of $\widetilde{V}$.
The right hand of \eqref{eqP<V3} is proved. The left hand is obvious by noting $Q(B)=1-Q(B^c)$ and $\widetilde{v}(B)=1-\widetilde{V}(B^c)$.
  The proof is completed. $\Box$

The second lemma is the Borel-Cantelli lemma for a countably sub-additive capacity.
\begin{lemma}\label{lemBCdirect} Let $V$ be a   countably sub-additive capacity and $\sum_{n=1}^{\infty}V(A_n)<\infty$. Then
$$ V(A_n\; i.o.)=0, \;\; \text{ where } \{ A_n\; i.o.\}=\bigcap_{n=1}^{\infty}\bigcup_{i=n}^{\infty}A_i. $$
\end{lemma}
{\bf Proof}. Easy and omitted.  $\Box$

The third  lemma is the converse part of  Borel-Cantelli lemma under $\widetilde{\Capc}^{\ast}$ or $\widetilde{V}$.
\begin{lemma}\label{lemBC}
 Let $\{X_n;n\ge 1\}$ be a sequence of independent random variables in the sub-linear expectation space $(\Omega,\mathscr{H},\Sbep)$ for which each $X_n$ is tight,
 $\{\tilde X_n;n\ge 1\}$ be  its copy  on   $(\widetilde{\Omega},\widetilde{\mathscr{H}},\widetilde{\mathbb E})$ as defined in Proposition \ref{lem3.4}.
 Suppose $\bm \xi_k=(X_{n_{k-1}+1},\ldots,X_{n_k})$, $1=n_0<n_1<\ldots$, $f_{k,j}\in C_{l,Lip}(\mathbb R^{n_k-n_{k-1}})$ and $\sum\limits_{k=1}^{\infty}\Capc( f_{k,j}(\bm\xi_k)\ge 1+\epsilon_{k,j} )=\infty$, $j=1,2,\ldots$.
   Then on  the space $(\widetilde{\Omega},\widetilde{\mathscr{H}},\widetilde{\mathbb E})$,
$$ \widetilde{\Capc}\left (  A\right)=\widetilde{\Capc}^{\ast}\left (A\right)=\widetilde{V}\left (  A\right)=1,\;\; A=\bigcap_{j=1}^{\infty}\big\{f_{k,j}(\tilde{\bm\xi}_k)\ge 1 \;\; i.o.\big\}, $$
 where $\tilde{\bm \xi}_k=(\tilde X_{n_{k-1}+1},\ldots,\tilde X_{n_k})$.
\end{lemma}

{\bf Proof}. Let $g_{k,j}\in C_{b,Lip}(\mathbb R)$ such that $I\{x\ge 1\}\ge g_{k,j}(x)\ge I\{x\ge 1+\epsilon_{k,j}\}$. Then
$$ \sum_{k=1}^{\infty} \Sbep\left[g_{k,j}\big(f_{k,j}(\bm \xi_k)\big)\right]=\infty, \; \; j=1,2,\ldots. $$
 By the expression \eqref{linearexpression}, for each pair of $k$ and $j$ there exists $\theta_{k,j}\in\Theta$ such that
$$ E_{\theta_{k,j}}\left[g_{k,j}\big(f_{j,k}(\bm \xi_k)\big)\right]=\Sbep\left[g_{k,j}\big(f_{j,k}(\bm \xi_k)\big)\right]. $$
Define the linear operator $E_k$ by
 $$E_k=\sum_{j=1}^{\infty} 2^{-j}E_{\theta_{k,j}}. $$
Then $E_k\le \Sbep$. By Proportion \ref{lem3.4} (d),  there exists a probability measure $Q$ on $\sigma(\tilde X_1,\tilde X_2,\ldots)$ such that $\{\bm \xi_k;k\ge 1\}$ is a sequence of independent random variables under $Q$, and \eqref{eqP<V1}-\eqref{eqP<V3} hold. By \eqref{eqP<V1},
\begin{align*}
 & \sum_{k=1}^{\infty} Q(f_{k,j}(\tilde{\bm \xi}_k)\ge 1)\ge  \sum_{n=1}^{\infty}Q\left[g_{k,j}\big(f_{k,j}(\tilde{\bm \xi}_k)\big)\right]=\sum_{k=1}^{\infty}E_k\left[g_{k,j}\big(f_{k,j}(\bm \xi_k)\big)\right]\\
& \;\; \ge  \frac{1}{2^j}\sum_{k=1}^{\infty}E_{\theta_{k,j}}\left[g_{k,j}\big(f_{k,j}(\bm \xi_k)\big)\right]
=\frac{1}{2^j}\sum_{k=1}^{\infty}\Sbep\left[g_{k,j}\big(f_{k,j}(\bm \xi_k)\big)\right]=\infty.
\end{align*}
So, by the Borell-Cantelli lemma for a  probability measure,
$$
Q\left ( f_{k,j}(\bm \xi_k)\ge 1 \;\; i.o. \right)=1.
$$
It follows that
$$ Q\left ( \bigcap_{j=1}^{\infty}\big\{f_{k,j}(\tilde{\bm \xi}_k)\ge 1 \;\; i.o.\big\}\right)=1. $$
By  \eqref{eqP<V3}, it follows that
$$ \widetilde{\Capc}\left (  A\right)\ge \widetilde{\Capc}^{\ast}\left (A\right)\ge \widetilde{V}\left (  A\right)=1. $$
The proof is now completed.
 $\Box$.

The next lemma tells us that the converse part of the Borel-Cantelli lemma remains valid in the original  sub-linear expectation space $(\Omega,\mathscr{H},\Sbep)$ under certain conditions.
\begin{lemma}\label{lemBC2} Let $(\Omega,\mathscr{H},\Sbep)$ be a sub-linear expectation  space with  a capacity $V$ having  the property (\ref{eq1.3}), and $v=1-V$. Suppose  one   of the following conditions is satisfied.
\begin{description}
      \item[\rm (a)]      The sub-linear expectation $\Sbep$ satisfies
$$ \Sbep[X]= \max_{P\in \mathscr{P}}P[X], \;\; X\in \mathscr{H}_b,$$
where $\mathscr{H}_b=\{f\in \mathscr{H}; f \text{ is bounded}\}$, $\mathscr{P}$ is a  countable-dimensionally weakly compact family of probability measures on $(\Omega,\sigma(\mathscr{H}))$ in sense that for any  $Y_1,Y_2,\ldots \in \mathscr{H}_b$ and any  sequence $\{P_n\}\subset \mathscr{P}$ there is a subsequence $\{n_k\}$ and a probability measure $P\in \mathscr{P}$ for which
\begin{equation}\label{eqcompact1} \lim_{k\to \infty} P_{n_k}[\varphi(Y_1,\ldots,Y_d)]= P[\varphi(Y_1,\ldots,Y_d)],\; \varphi\in C_{b,Lip}(\mathbb R^d), d\ge 1.
\end{equation}
 \item[\rm (b)]  $\Sbep$ on $\mathscr{H}_b$ is regular in sense that $\Sbep[X_n]\downarrow 0$  for any   elements  $\mathscr{H}_b\ni X_n\downarrow0$. Let $\mathscr{P}$ be the family of all  probability measures on $(\Omega, \sigma(\mathscr{H}))$ for which
$$P[f]\le \Sbep[f], \;\; f\in \mathscr{H}_b. $$
  \item[\rm (c)]     $\Omega$ is a complete separable metric space, each element $X(\omega)$ in $\mathscr{H}$ is a continuous function on $\Omega$.  The capacity $V$ with the property (\ref{eq1.3})   is tight in sense that, for any $\epsilon>0$  there is a compact set $K\subset \Omega$ such that $V(K^c)<\epsilon$.  Let $\mathscr{P}$ be defined as in (b).
  \item[\rm (d)]  $\Omega$ is a complete separable metric space, each element $X(\omega)$ in $\mathscr{H}$ is a continuous function on $\Omega$. The sub-linear expectation $\Sbep$ is defined by
$$ \Sbep[X]= \max_{P\in \mathscr{P}}P[X], $$
where $\mathscr{P}$ is a weakly compact family of probability measures on $(\Omega,\mathscr{B}(\Omega))$.
\end{description}
Denote $\Capc^{\mathscr{P}}(A)=\max_{P\in \mathscr{P}}P(A)$, $A\in \sigma(\mathscr{H}). $
 Let $\{X_n;n\ge 1\}$ be a sequence of independent random variables in  $(\Omega,\mathscr{H},\Sbep)$.
 \begin{description}
   \item[\rm (i)] If $\sum_{n=1}^{\infty}v(X_n<1)<\infty$,  then for $\Capc=\Capc^{\mathscr{P}}$,  $\mathbb C^{\ast}$,  $\outCapc$ or $\upCapc$,
   \begin{equation}\label{eqlemBC2.1}\Capc\left(\bigcup_{m=1}^{\infty} \bigcap_{i=m}^n \{X_i\ge 1\}\right)=1\;\; i.e., \; \cCapc\left( X_i<1 \; i.o.\right)=0.
   \end{equation}
   \item[\rm (ii)]   If $\sum_{n=1}^{\infty}V(X_n\ge 1)=\infty$, then for $\Capc=\Capc^{\mathscr{P}}$,  $\mathbb C^{\ast}$,  $\outCapc$ or $\upCapc$,
   \begin{equation}\label{eqlemBC2.2} \Capc\left( X_n\ge 1\;\; i.o. \right)=1.
   \end{equation}
   More generally, suppose $\{\bm X_n;n\ge 1\}$ is a sequence of independent random vectors  in $(\Omega,\mathscr{H},\Sbep)$, where $\bm X_n$ is $d_n$-dimensional,  $f_{n,j}\in C_{l,lip}(\mathbb R^{d_n})$ and $\sum_{n=1}^{\infty}V(f_{n,j}(\bm X_n)\ge 1)=\infty$, $j=1,2,\ldots$,  then for $\Capc=\Capc^{\mathscr{P}}$,  $\mathbb C^{\ast}$,  $\outCapc$ or $\upCapc$,
    \begin{equation}\label{eqlemBC2.3} \Capc\left( \bigcap_{j=1}^{\infty}\big\{f_{n,j}(\bm X_n)\ge 1\;\; i.o.\big\} \right)=1.
   \end{equation}
\item[\rm (iii)] Suppose $\{\bm X_n;n\ge 1\}$ is a sequence of independent random vectors  in $(\Omega,\mathscr{H},\Sbep)$, where $\bm X_n$ is $d_n$-dimensional. If $F_n$ is a $d_n$-dimensional close set with $\sum_{n=1}^{\infty} v(\bm X_n \not\in F_n)<\infty$, then for $\Capc=\Capc^{\mathscr{P}}$,  $\mathbb C^{\ast}$,  $\outCapc$ or $\upCapc$,
          $$ \cCapc\left( \bm X_n\not\in F_n\; i.o.\right)=0; $$
          If $F_{n,j}$ are $d_n$-dimensional  close sets with $\sum_{n=1}^{\infty} V(\bm X_n\in F_{n,j})=\infty$, $j=1,2,\ldots$, then for $\Capc=\Capc^{\mathscr{P}}$,  $\mathbb C^{\ast}$,  $\outCapc$ or $\upCapc$,
          $$ \Capc\left(\bigcap_{j=1}^{\infty}\big\{\bm X_n\in F_{n,j}\;\; i.o.\big\}\right)=1. $$
 \end{description}
\end{lemma}

{\bf Proof.} (i) and (ii) are special cases of (iii). But, to prove the general case (iii), we need to show the two special cases first.   Without loss of generality, we can assume $0\le X_n\le 2$, for otherwise, we can replace it by $0\vee X_n\wedge 2$. Write $\bm X=(X_1,X_2,\ldots)$. Suppose (a) is satisfied. Consider the family $\mathscr{P}$ on $\sigma(\bm X)$. We denote it by  $\mathscr{P}\big|_{\sigma(\bm X)}$. Note $|X_n|\le 2$, $n=1,2,\ldots$, and the set $K=\bigotimes_{i=1}^{\infty}[-2,2]$ is a compact set on $\mathbb R^{\infty}$. So,  $\mathscr{P}\big|_{\sigma(\bm X)}$ is a tight and so a relatively weakly compact family, i.e., the family $\{\overline{P}:\overline{P}(A)=P(\bm X\in A), A\in \mathscr{B}(\mathbb R^{\infty}), P\in \mathscr{P}\}$ is  a relatively weakly compact family of probability measures on the metric space $\mathbb R^{\infty}$. Next, we show that  $\mathscr{P}\big|_{\sigma(\bm X)}$ is close. Suppose $P_n$ is  weakly convergent on  $\sigma(\bm X)$.
Then there exists a probability on $\mathbb R^{\infty}$ such that $\overline{P}_n\implies Q$, i.e.,
\begin{equation}\label{eqlimit-P} Q[f]=\lim_{n\to \infty}P_n[f(\bm X], f\in C_b(\mathbb R^{\infty}).
\end{equation}
 It is needed to show that there exists  a probability measure $P\in \mathscr{P}$ satisfying $Q(A)=P\big(\bm X\in A\big)$ for  $A\in \mathscr{B}(\mathbb R^{\infty})$.
 By the conditions assumed, for the sequence $\{P_n\}$ there exists a subsequence $\{n_k\}$ and a probability measure $P\in \mathscr{P}$ such that \eqref{eqcompact1} holds.  Hence
 $$ Q[f]=P[f(X_1,\ldots,X_d)], \;\; \forall f\in C_{b,lip}(\mathbb R^d),\; d\ge 1. $$
So, $Q(\{\bm x: (x_1,\ldots, x_d)\in A\})=P\big((X_1,\ldots, X_d)\in A\big)$ for all $A\in \mathscr{B}(\mathbb R^d)$, which implies  $Q(A)=P\big(\bm X\in A\big)$ for all $A\in \mathscr{B}(\mathbb R^{\infty})$.
We conclude that $\mathscr{P}\big|_{\sigma(\bm X)}$ is close and so weakly compact. Denote $\widetilde{V}(A)=\Capc^{\mathscr{P}}(\bm X\in A)$. By Lemma 6.1.12 of Peng (2019), for any sequence of closed sets $F_n\downarrow F$, we have $\widetilde{V}( F_n)\downarrow \widetilde{V}(F)$.

Now, we consider (i). By the independence, we have for any $\delta_i>0$, and $\Capc=\Capc^{\mathscr{P}}$,  $\mathbb C^{\ast}$,  $\outCapc$ or $\upCapc$,
$$ \Capc\left(\bigcap_{i=m}^n \{X_i\ge 1-\delta_i\}\right)\ge \prod_{i=m}^n V(X_i\ge 1). $$
In fact, we can choose a Lipschitz function $f_i$ such that $I\{x\ge 1-\delta_i\}\ge f_i(x)\ge I\{x\ge 1\}$. Then
$$ \Capc\left(\bigcap_{i=m}^n \{X_i\ge 1-\delta_i\}\right)\ge \Sbep[\prod_{i=m}^n f_i(X_i)]=\prod_{i=m}^n \Sbep[f_i(X_i)]\ge \prod_{i=m}^n V(X_i\ge 1). $$
Let $\epsilon_i=v(X_i<1)$ and   choose $\delta_i=1/l$. Then
$$ \Capc^{\mathscr{P}}\left(\bigcap_{i=m}^n \{X_i\ge 1-1/l\}\right)\ge \prod_{i=m}^{\infty} V(X_i\ge 1) =\prod_{i=m}^{\infty}(1-\epsilon_i). $$
Note $\left\{\bm x: \bigcap_{i=m}^n \{x_i\ge 1-\delta_i\}\right\}$ is a close set of $\bm x$ on $\mathbb R^{\infty}$. It follows that
$$ \Capc^{\mathscr{P}}\left(\bigcap_{i=m}^n \{X_i\ge 1-1/l\}\right)\searrow_l \Capc^{\mathscr{P}}\left(\bigcap_{i=m}^n \{X_i\ge 1\}\right)\searrow_n \Capc^{\mathscr{P}}\left(\bigcap_{i=m}^{\infty} \{X_i\ge 1\}\right). $$
It follows that
$$ \Capc^{\mathscr{P}}\left(\bigcap_{i=m}^{\infty} \{X_i\ge 1\}\right)\ge \prod_{i=m}^{\infty} V(X_i\ge 1) =\prod_{i=m}^{\infty}(1-\epsilon_i)\to 1,\; m\to \infty $$
due the fact that $\sum_{i=1}^{\infty}\epsilon_i<\infty$. Hence \eqref{eqlemBC2.1} is proved.

Consider (ii). Write $\epsilon_i=V(X_i\ge 1)$. Now, for $\cCapc=1-\Capc^{\mathscr{P}}$,  $1-\mathbb C^{\ast}$,  $1-\outCapc$ or $1-\upCapc$, we have
$$ \cCapc\left(\bigcap_{i=m}^n \{X_i< 1-1/l\}\right)\le \cSbep\Big[ \prod_{i=m}^n\big(1-f_i(X_i)\big)\Big] =  \prod_{i=m}^n\cSbep\big[1-f_i(X_i)\big]\le  \prod_{i=m}^n v(X_i< 1). $$
That is
$$ \Capc\left(\bigcup_{i=m}^n \{X_i\ge  1-1/l\}\right)\ge 1-\prod_{i=m}^n \big(1-V(X_i\ge 1)\big)\ge 1-\exp\{-\sum_{i=m}^n \epsilon_i\}. $$
Note that $\bigcup_{i=m}^n \{x_i\ge  1-1/l\}$ is a close set of $\bm x$. It follows that
$$ \Capc^{\mathscr{P}}\left(\bigcup_{i=m}^n \{X_i\ge  1-1/l\}\right)\searrow \Capc^{\mathscr{P}}\left(\bigcup_{i=m}^n \{X_i\ge  1\}\right)\; \text{ as } l\to \infty. $$
Hence for each $m$,
\begin{equation}\label{eqprooflemBC2.3} \Capc^{\mathscr{P}}\left(\bigcup_{i=m}^n \{X_i\ge  1\}\right)\ge 1-\exp\{-\sum_{i=m}^n \epsilon_i\}\to 1, \; n\to \infty,
\end{equation}
due to the fact that $\sum_{i=1}^{\infty}\epsilon_i=\infty$.
Let $\delta_k=2^{-k}$. We can choose a sequence $n_k\nearrow \infty$ such that
$$ \Capc^{\mathscr{P}}\left( \max_{n_k+1\le i\le n_{k+1}}X_i\ge  1\right)=\Capc^{\mathscr{P}}\left(\bigcup_{i=n_k+1}^{n_{k+1}} \{X_i\ge  1\}\right)\ge 1-\delta_k. $$
Let $Z_k=\max_{n_k+1\le i\le n_{k+1}}X_i$. Then $\{Z_k;k\ge 1\}$ are independent under $\Sbep$.  By (i),
$$ \Capc^{\mathscr{P}}\left(\bigcup_{l=1}^{\infty}\bigcap_{k=l}^{\infty}\{Z_k\ge 1\}\right)=1. $$
Note $\bigcup_{l=1}^{\infty}\bigcap_{k=l}^{\infty}\{Z_k\ge 1\}\subset \{X_n\ge 1\; i.o.\}$. \eqref{eqlemBC2.2} holds.

Now, we consider the general case. Without loss of generality, assume $0\le f_{n,j}(\bm X_n)\le 2$.  Similar to \eqref{eqprooflemBC2.3}, for each $m$ and $j$ we have
$$ \Capc^{\mathscr{P}}\left(\bigcup_{i=m}^n \{f_{i,j}(\bm X_i)\ge  1\}\right)\ge 1-\exp\{-\sum_{i=m}^n V(f_{i,j}(\bm X_i)\ge  1)\}\to 1, \; n\to \infty.
$$
Let $\delta_k=2^{-k}$. We choose the sequence $1=n_{0,0}<n_{1,1}<n_{2,1}<n_{2,2}<\ldots<n_{k,1}<\ldots<n_{k,k}<n_{k+1,1}<\ldots$ such that
$$ \Capc^{\mathscr{P}}\left(\bigcup_{i=n_{k,j-1}+1}^{n_{k,j}} \{f_{i,j}(\bm X_i)\ge  1\}\right)\ge 1-\delta_{k+j},\;\; j\le k, k\ge 1,  $$
where $n_{k,0}=n_{k-1,k-1}$. Let $Z_{k,j}=\max_{n_{k,j-1}+1\le i\le n_{k,j}}f_{i,j}(\bm X_i)$. Then the random variables $Z_{1,1},Z_{2,1},Z_{2,2},\ldots,Z_{k,1},\ldots,Z_{k,k},Z_{k+1,1},\ldots$ are independent under $\Sbep$ with
$$ \cCapc^{\mathscr{P}}(Z_{k,j}<1)< \delta_{k+j}. $$
Note $\sum_{k=1}^{\infty}\sum_{j=1}^k\delta_{k+j}<\infty$. By (i), we have
$$ \Capc^{\mathscr{P}}\left(\bigcup_{l=1}^{\infty}\bigcap_{k=l}^{\infty}\bigcap_{j=1}^k\{Z_{k,j}\ge 1\}\right)=1. $$
On the event $\bigcup_{l=1}^{\infty}\bigcap_{k=l}^{\infty}\bigcap_{j=1}^k\{Z_{k,j}\ge 1\}$, there exists a $l_0$ such that $Z_{k,j}\ge 1$ for all $k\ge l_0$ and $1\le j\le k$. For each fixed $j$, when $k\ge j\vee l_0$ we have $Z_{k,j}\ge 1$, and hence $\{f_{n,j}(\bm X_n)\ge  1 \;\; i.o\}$ occurs.
It follows that
$$ \bigcup_{l=1}^{\infty}\bigcap_{k=l}^{\infty}\bigcap_{j=1}^k\{Z_{k,j}\ge 1\}\subset \bigcap_{j=1}^{\infty}\{f_{n,j}(\bm X_n)\ge  1 \;\; i.o\}. $$
\eqref{eqlemBC2.3} holds.

(iii) Denote $d(\bm x,F)=\inf\{\|\bm y-\bm x\|: \bm y\in F\}$. Then $d(\bm x,F)$ is a Lipschitz function of $\bm x$. If $F_{n,j}$ is a close set, then
$$ \bm X_n\in F_{n,j}\Longleftrightarrow d(\bm X_n, F_{n,j})=0\Longleftrightarrow f_{n,j}(\bm X_n)=:1-1\wedge d(\bm X_n, F_{n,j})\ge 1. $$
The results follows from (i) and (ii) immediately.

When the condition (b) is satisfied, it is sufficient to show that the family $\mathscr{P}$ satisfies the assumption in (a). Note the  expression (\ref{linearexpression}).
Consider the linear expectation $E_{\theta}$ on $\mathscr{H}_b$. If $\mathscr{H}_b\ni f_n\downarrow 0$, then $0\le E_{\theta}[f_n]\le \Sbep[f_n]\to 0$.
Hence, similar to Lemma 1.3.5 and Lemma 6.2.2 of Peng (2019), by Daniell-Stone's theorem, there is a unique probability $P_{\theta}$ on $\sigma(\mathscr{H}_b)=\sigma(\mathscr{H})$ such that
 $$ E_{\theta}[f]=   P_{\theta}[f]\le \Sbep[f],  \; f\in \mathscr{H}_b. $$
 Hence
 $$\Sbep[f]=\sup_{\theta\in \Theta}E_{\theta}[f]=  \sup_{\theta\in \Theta} P_{\theta}[f],  f\in \mathscr{H}_b. $$
Recall that $\mathscr{P}$ is the family of all probability measures $P$ on $\sigma(\mathscr{H})$ which satisfies $P[f]\le \Sbep[f]$ for all $f\in \mathscr{H}_b$.
We have
 $$\Sbep[f]=  \sup_{\theta\in \Theta} P_{\theta}[f]\le\sup_{P\in \mathscr{P}}P[f]\le \Sbep[f],\;\;  f\in \mathscr{H}_b. $$
Suppose $Y_1,Y_2,\ldots\in \mathscr{H}_b$ with $|Y_i|\le C_i$. Write $\bm Y=(Y_1,Y_2,\ldots)$ and $K=\bigotimes_{i=1}^{\infty}[-C_i,C_i]$. Then $P(\bm Y\in K^c)=0$ and $K$ is a compact set on the space $\mathbb R^{\infty}$. It follows that $\mathscr{P}\big|_{\sigma(\bm Y)}$ is tight and so is relatively weakly compact. Hence, for any sequence $\{P_n\}\subset \mathscr{P}$, there exists a subsequence $n_k\nearrow \infty$ such that
$$E[f(\bm Y)]=\lim_{k\to \infty}P_{n_k}[f(\bm Y)], \; f\in C_b(\mathbb R^{\infty})$$
is well-defined. It is obvious that $E$ is a linear expectation on  $\{\varphi(\bm Y):\varphi\in C_b(\mathbb R^{\infty})\}$. Consider $E$ on $\mathscr{L}=\{\varphi(Y_1,\ldots,Y_d): \varphi\in C_{b,Lip}(\mathbb R^d), d\ge 1\}$. It is obvious that
$$ E[\varphi(Y_1,\ldots,Y_d)]= \lim_{k\to \infty}P_{n_k}[\varphi(Y_1,\ldots,Y_d)]\le \Sbep[[f(Y_1,\ldots,Y_d)],\; \varphi\in C_{b,Lip}(\mathbb R^d). $$
So,  by the Hahn-Banach theorem, there exists  a finite additive linear expectation  $E^e$ defined on $\mathscr{H}$ such that, $E^e=E$ on $\mathscr{L}$ and, $E^e\le \Sbep$ on $\mathscr{H}$. For $E^e$, by the regularity, as shown before there is probability measure $P^e$ on $\sigma(\mathscr{H})$ such that $P^e[f]=E^e[f]$ for all $f\in \mathscr{H}_b\supset \mathscr{L}$. Hence $P^e\in \mathscr{P}$ and
$$ \lim_{k\to \infty}P_{n_k}[\varphi(Y_1,\ldots,Y_d)]=E[\varphi(Y_1,\ldots,Y_d)]=P^e[\varphi(Y_1,\ldots,Y_d)], \;\; \varphi\in C_{b,Lip}(\mathbb R^d), d\ge 1. $$
It follows that $\mathscr{P}$ satisfies the assumption in (a).

For (c), it can be shown  that $\Sbep$ is regular on $\mathscr{H}_b$ and so the condition (b) is satisfied. Also, (d) is a special case of (a). The proof is completed.
$\Box$

\bigskip

The  rest   three lemmas give the estimators of the tail capacities of maximum partial sums of independent random variables. Lemma \ref{moment_v} below is a kind of Kolmogorov's maximal inequality  under  $\lowCapc$.
\begin{lemma}\label{moment_v}
Let    $\{\bm Z_{n,k}; k=1,\ldots, k_n\}$ be an array of independent random vectors taking values in $\mathbb R^d$    such that $\Sbep[|\bm Z_{n,k}|^2]<\infty$, $k=1,\ldots, k_n$, here $|\cdot|$ is the Euclidean norm. Then for any $\bm\mu_{n,k}\in \widetilde{\mathbb M}[\bm Z_{n,k}]=:\big\{E[\bm Z_{n,k}]: E\in \mathscr{E}\big\}$ where $\mathscr{E}$  is defined as \eqref{linearexpression2.1},   $k=1,\ldots, k_n$,
\begin{equation*}
\lowCapc\left(\max_{m\le k_n}\big|\sum_{k=1}^{m} (\bm Z_{n,k}-\bm \mu_{n,k}])\big| \geq x\right)\leq 2x^{-2}\sum_{k=1}^{k_n}\left(\Sbep[|\bm Z_{n,k}|^2]-|\bm \mu_{n,k}|^2\right) ,\ \ \forall x>0.
\end{equation*}
\end{lemma}
{\bf Proof}. For each $k$ there   exists  $E_k\in \mathscr{E}$ such that $\bm \mu_{n,k}=E_k[\bm Z_{n,k}]$.    $Z_k$ is a finite additive linear expectation on $\mathscr{H}_b=\{f\in \mathscr{H}; f \text{ is bounded}\}$ with $E\le \vSbep=\Sbep$. Note that each $\bm Z_{n,k}$ is tight by the fact $\Sbep[|\bm Z_{n,k}|^2]<\infty$.
By Proposition \ref{lem3.4},  $\{\bm Z_{n,k}; k=1,\ldots, k_n\}$ has a copy $\{\tilde{\bm Z}_{n,k}; k=1,\ldots, k_n\}$  on a new sub-linear expectation space $(\widetilde{\Omega},\widetilde{\mathscr{H}}, \widetilde{\mathbb E})$  with a probability measure $Q$ on $\sigma(\tilde{\bm Z}_{n,1},\ldots,\tilde{\bm Z}_{n,k_n})$ such that $\{\tilde{\bm Z}_{n,1},\ldots, \tilde{\bm Z}_{n,k_n}\}$ are independent random vectors under $Q$,
\begin{equation} \label{eqproofmoment-v.1}
Q\big[\varphi(\tilde{\bm Z}_{n,k})\big]=E_k\left[\varphi(\bm Z_{n,k})\right]\; \text{ for all } \varphi\in C_{b,Lip}(\mathbb R^d),
\end{equation}
\begin{equation}  \label{eqproofmoment-v.2}
Q\left[\varphi(\tilde{\bm Z}_{n,1},\ldots,\tilde{\bm Z}_{n,k_n})\right]\le \Sbep\left[\varphi(\bm Z_{n,1},\ldots,\bm Z_{n,k_n})\right]\; \text{ for all } \varphi\in C_{b,Lip}(\mathbb R^{d\times k_n})
\end{equation}
and
\begin{equation}  \label{eqproofmoment-v.3}
\widetilde{\cCapc}(B)\le Q(B)\le \widetilde{\Capc}(B) \; \text{ for all } B\in\sigma(\tilde{\bm Z}_{n,1},\ldots,\tilde{\bm Z}_{n,k_n}).
\end{equation}
Note $E_k|Z_{n,k,i}-(-c)\vee Z_{n,k,i}\wedge c|\le \Sbep[|Z_{n,k,i}|-c)^+]\to 0$ as $c\to \infty$ by $\Sbep[|Z_{n,k}|^2]<\infty$. Then
$$Q[\tilde Z_{n,k,i}]=\lim\limits_{c\to \infty}Q[(-c)\vee \tilde Z_{n,k,i}\wedge c]
=\lim\limits_{c\to \infty}E_k[(-c)\vee Z_{n,k,i}\wedge c]=E_k [  Z_{n,k,i} ] =\mu_{n,k,i} $$
by \eqref{eqproofmoment-v.1}, and
$$ Q[|\tilde{\bm Z}_{n,k}|^2]=\lim\limits_{c\to \infty}Q[|\tilde{\bm Z}_{n,k}|^2\wedge c]\le \lim\limits_{c\to \infty}\Sbep[|\bm Z_{n,k}|^2\wedge c]
\le \Sbep[|\bm Z_{n,k}|^2]. $$
by \eqref{eqproofmoment-v.2}. Let $Y=\max\limits_{m\le k_n}\big|\sum\limits_{k=1}^{m}(\tilde{\bm Z}_{n,k}-\bm \mu_{n,k})\big|$.
By \eqref{eqproofmoment-v.3}   and the Kolmogorov inequality for independent random variables in a probability space, we have
\begin{align*}
  & \widetilde{\cCapc}\left(Y\ge x\right)\le \widetilde{v}\left(Y\ge x\right)\le Q\left(Y\ge x\right)
 \le     2x^{-2}\sum_{k=1}^{k_n}Q[|\tilde{\bm Z}_{n,k}-Q [\tilde{\bm Z}_{n,k}]|^2]\\
& =   2x^{-2}\sum_{k=1}^{k_n}(Q[|\tilde{\bm Z}_{n,k}|^2]-\big|Q [\tilde{\bm Z}_{n,k}]\big|^2)
\le 2 x^{-2}\sum_{k=1}^{k_n}\left(\Sbep[|\bm Z_{n,k}|^2]-|\bm \mu_{n,k}|^2\right).
\end{align*}
  By \eqref{eqV-V} and noting $\max\limits_{m\le k_n}\big|\sum\limits_{k=1}^{m}(\tilde{\bm Z}_{n,k}-\bm \mu_{n,k})\big|\overset{d}=\max\limits_{m\le k_n}\big|\sum\limits_{k=1}^{m}(\bm Z_{n,k}-\bm \mu_{n,k})\big|$, we have
\begin{align*}
&\lowCapc\left( \max\limits_{m\le k_n}\big|\sum\limits_{k=1}^{m}(\bm Z_{n,k}-\bm \mu_{n,k})\big|\ge x\right)\\
\le &  \widetilde{\cCapc}\left(Y\ge y \right)\le 2 y^{-2}\sum_{k=1}^{k_n}\left(\Sbep[|\bm Z_{n,k}|^2]-|\bm \mu_{n,k}|^2\right), \;\; 0<y<x.
\end{align*}
 The proof is completed. $\Box$

The following lemma is on the exponential inequality under $\upCapc$ whose proof is similar to that of Theorem 4.5 of Zhang (2021a).
\begin{lemma}\label{lem4.1}   Let $\{Z_{n,k}; k=1,\ldots, k_n\}$ be an array of independent random variables under $\Sbep$ such that $\Sbep[Z_{n,k}]\le 0$ and $\Sbep[Z_{n,k}^2]<\infty$, $k=1,\ldots, k_n$.     Then for all $x,y>0$
 \begin{align}\label{eqlem4.1.1}
 &\upCapc\left(\max_{m\le k_n}\sum_{k=1}^{m} Z_{n,k}\ge x\right)\nonumber\\
 \le &\upCapc\left(\max_{k\le k_n} Z_{n,k}\ge y\right) + \exp\left\{\frac{x}{y}-\frac{x}{y}\Big(\frac{B_n^2}{xy}+1\Big)\ln\Big(1+\frac{xy}{B_n^2}\Big)\right\},
 \end{align}
 where $B_n^2=\sum_{k=1}^{k_n}\Sbep[Z_{n,k}^2]$. In particular, by letting $y=x$, we have Kolmogorov's maximal inequality under $\upCapc$ as follows:
 \begin{equation}\label{eqlem4.1.2}
 \upCapc\left(\max_{m\le k_n}\sum_{k=1}^{m} Z_{n,k}\ge x\right)\le (e+1)\frac{B_n^2}{x^2},\ \ \forall x>0.
 \end{equation}
 \end{lemma}

     The last  lemma    on  the L\'evy maximal inequality is Lemma 2.1 of Zhang (2020).

\begin{lemma}\label{LevyIneq} Let $X_1,\cdots, X_n$ be independent random variables in  a sub-linear expectation space $(\Omega, \mathscr{H}, \Sbep)$, $S_k=\sum_{i=1}^k X_i$, and $0<\alpha<1$ be a real number.
If there exist real constants $\beta_{n, k}$  such that
$$ \upCapc\left(|S_k-S_n|\ge \beta_{n,k}+\epsilon\right)\le \alpha, \text{ for all } \epsilon>0 \text{ and } k=1,\cdots ,n, $$
then
\begin{equation}\label{eqLIQ2}   (1-\alpha) \upCapc\left(\max_{k\le n}(|S_k| -\beta_{n,k})> x+\epsilon\right)\le   \upCapc\left(|S_n|>x\right), \text{ for all }x>0, \epsilon>0.
\end{equation}
\end{lemma}

 \section{The law of large numbers} \label{sectMain}
\setcounter{equation}{0}
Our first theorem gives the sufficient and necessary conditions for the strong law of large numbers without any assumption on the continuity of the capacities. Let $\{X_n;n\ge 1\}$ be a sequence of i.i.d. random variables in a sub-linear expectation space $(\Omega,\mathscr{H},\Sbep)$. Denote  $S_n=\sum_{i=1}^n X_i$.
\begin{theorem}\label{thLLN1}
\begin{description}
  \item[\rm (a)]  If
  \begin{equation}\label{eqthLLN1.1} C_{\upCapc}(|X_1|)<\infty,
  \end{equation}
  then
  \begin{equation}\label{eqthLLN1.2}
  \outCapc\left(\liminf_{n\to\infty}\frac{S_n}{n}<\vcSbep[X_1]\;\text{ or }\; \limsup_{n\to\infty} \frac{S_n}{n}>\vSbep[X_1]\right)=0.
  \end{equation}
  Further, if the space $(\Omega,\mathscr{H},\Sbep)$ satisfies one of the conditions (a)-(d) in Lemma \ref{lemBC2}, then for $\Capc=\Capc^{\mathscr{P}}$, $\mathbb C^{\ast}$, $\outCapc$ or $\upCapc$,
   \begin{align}\label{eqthLLN1.3}
   &\Capc\left(\liminf_{n\to\infty}\frac{ S_n}{n}=\vcSbep[X_1]\;\text{ and }\; \limsup_{n\to\infty} \frac{ S_n}{n}=\vSbep[X_1]\right)=1,
  \end{align}
  \begin{equation}\label{eqthLLN1.4}
  \Capc\left(C\left\{\frac{ S_n}{n}\right\} =[\vcSbep[X_1],\vSbep[X_1]]\right)=1,
  \end{equation}
  where $C\{x_n\}$ denotes the cluster set of a sequence of $\{x_n\}$ in $\mathbb R$.
  \item[\rm (b)] Suppose the space $(\Omega,\mathscr{H},\Sbep)$ satisfies one of the conditions (a)-(d) in Lemma \ref{lemBC2}. If   for $\Capc=\Capc^{\mathscr{P}}$, $\mathbb C^{\ast}$,  $\outCapc$ or $\upCapc$,
  \begin{equation}\label{eqthLLN1.5}
  \Capc\left(\limsup_{n\to\infty}\frac{| S_n|}{n}=\infty\right)<1,
  \end{equation}
  then \eqref{eqthLLN1.1} holds.
\end{description}

\end{theorem}

\begin{remark} Theorem \ref{thLLN1} tells us that the sufficient and necessary condition for the strong law of large numbers is \eqref{eqthLLN1.1}. Under
\eqref{eqthLLN1.1}, $\vSbep[X_1]$ and $\vcSbep[X_1]$ are well-defined and finite. In Zhang (2016), \eqref{eqthLLN1.2} is proved under \eqref{eqthLLN1.1} and an extra condition that $\Sbep[(|X_1|-c)^+]\to 0$ as $c\to \infty$. Under this extra condition, we have $\Sbep[X_1]=\vSbep[X_1]$ and $\cSbep[X_1]=\vcSbep[X_1]$. For establishing \eqref{eqthLLN1.3} and \eqref{eqthLLN1.4} and (b), the continuity of $\outCapc$ is also assumed in Zhang (2016).
\end{remark}

  The following corollary gives an analogues of \eqref{eqKol1}.

\begin{corollary} \label{cor1} Suppose the space $(\Omega,\mathscr{H},\Sbep)$ satisfies one of the conditions (a)-(d) in Lemma \ref{lemBC2}.
\begin{description}
  \item[\rm (a)]  If \eqref{eqthLLN1.1} is satisfied, then for $\Capc=\Capc^{\mathscr{P}}$, $\mathbb C^{\ast}$,  $\outCapc$ or $\upCapc$,
  \begin{align}\label{eqcor1.1}
   \Capc\left(\lim_{n\to\infty}\frac{  S_n}{n}=b\right)
  =\begin{cases} 1, & \text{ when } b\in \big[\vcSbep[X_1], \vSbep[X_1]\big],\\
  0, & \text{ when } b\not\in \big[\vcSbep[X_1], \vSbep[X_1]\big].
  \end{cases}
  \end{align}
  \item[\rm (b)]   For $\Capc=\Capc^{\mathscr{P}}$, $\mathbb C^{\ast}$,  $\outCapc$ or $\upCapc$, there exists a  finite random variable $b$ such that
  \begin{equation}\label{eqcor1.2}
   \cCapc\left(\lim_{n\to\infty}\frac{  S_n}{n}=b\right)=1
  \end{equation}
  if and only if \eqref{eqthLLN1.1}, $\vcSbep[X_1]=\vSbep[X_1]$ and $ \cCapc(b=\vSbep[X_1])=1$.
\end{description}
\end{corollary}

The following theorem and corollary are Marcinkiewicz's type laws of large numbers which gives the rate of convergence of Kolmogorov's type law of large numbers.
\begin{theorem}\label{thLLN2} Let $1\le p<2$.    If
  \begin{equation}\label{eqthLLN2.1} C_{\upCapc}(|X_1|^p)<\infty,
  \end{equation}
 then
  \begin{equation}\label{eqthLLN2.2}
  \outCapc\left(\liminf_{n\to\infty}\frac{S_n-n\vcSbep[X_1]}{n^{1/p}}<0\;\text{ or }\; \limsup_{n\to\infty} \frac{S_n-n\vSbep[X_1]}{n^{1/p}}>0\right)=0.
  \end{equation}
 Further, if the space $(\Omega,\mathscr{H},\Sbep)$ satisfies one of the conditions (a)-(d) in Lemma \ref{lemBC2}, then for $\Capc=\Capc^{\mathscr{P}}$, $\mathbb C^{\ast}$, $\outCapc$ or $\upCapc$,
   \begin{align}\label{eqthLLN2.3}
 \Capc\left(\liminf_{n\to\infty}\frac{ S_n-n\vcSbep[X_1]}{n^{1/p}}=0\;\text{ and }\; \limsup_{n\to\infty} \frac{ S_n-n\vSbep[X_1]}{n^{1/p}}=0\right)=1.
  \end{align}
\end{theorem}

\begin{corollary}\label{cor2} Suppose the space $(\Omega,\mathscr{H},\Sbep)$ satisfies one of the conditions (a)-(d) in Lemma \ref{lemBC2}.
\begin{description}
  \item[\rm (a)]  If \eqref{eqthLLN2.1}  is satisfied, then for $\Capc=\Capc^{\mathscr{P}}$, $\mathbb C^{\ast}$,  $\outCapc$ or $\upCapc$,
  \begin{align}\label{eqcor2.1}
   \Capc\left(\lim_{n\to\infty}\frac{ S_n-nb}{n^{1/p}}=0\right)
 =   \begin{cases} 1, & \text{ when } b\in \big[\vcSbep[X_1], \vSbep[X_1]\big],\\
  0, & \text{ when } b\not\in \big[\vcSbep[X_1], \vSbep[X_1]\big].
  \end{cases}
  \end{align}
  \item[\rm (b)]  For $\Capc=\Capc^{\mathscr{P}}$, $\mathbb C^{\ast}$,  $\outCapc$ or $\upCapc$, there exists   finite random variables $b(\omega)$ and $c(\omega)$ such that
  \begin{equation}\label{eqcor2.2}
 \cCapc\left(\lim_{n\to\infty}\frac{ S_n-nb}{n^{1/p}}=c\right)=1
  \end{equation}
  if and only if $C_{\upCapc}(|X_1|^p)<\infty$, $\vcSbep[X_1]=\vSbep[X_1]$ and
  \begin{equation}\label{eqcor2.3}
  \begin{aligned}
  \cCapc(b+c=\vSbep[X_1])=1, & \text{ when } p=1, \\
  \cCapc(b=\vSbep[X_1], c=0)=1, & \text{ when } 1<p<2.
   \end{aligned}
   \end{equation}
\end{description}
\end{corollary}

\begin{remark} When the space $(\Omega,\mathscr{H},\Sbep)$ does not satisfy  the conditions (a)-(d) in Lemma \ref{lemBC2}. We may consider the   the copy $\{\tilde X_n;n\ge 1\}$. The sub-linear expectation space $(\widetilde{\Omega}, \widetilde{\mathscr{H}},\widetilde{\mathbb E})$ satisfies the condition (d) in Lemma \ref{lemBC2} when each $X_n$ is tight, by Proposition \ref{lem3.4}, and so, Theorems \ref{thLLN1} and \ref{thLLN2} and Corollary \ref{cor1} and \ref{cor2} remain true for $\{\tilde X_n;n\ge 1\}$.
\end{remark}
 Zhang (2020) studied the the convergence of the infinite series $\sum_{n=1}^{\infty}X_n$ of a sequence of independent random variables. But, when the strong convergence is considered, the   capacity $\upCapc$ is assumed to be continuous. The next  theorem gives the equivalence among  various kinds of the convergence  of the infinite series $\sum_{n=1}^{\infty}X_n$ without the assumption   on the continuity of the capacities.

\begin{theorem}\label{th3}Let $\{X_n;n\ge 1\}$ be a sequence of independent random variables in a sub-linear expectation space $(\Omega,\mathscr{H},\Sbep)$ with a capacity   satisfying \eqref{eq1.4}, and  $\{\tilde X_n;n\ge 1\}$ be its copy on $(\widetilde{\Omega},\widetilde{\mathscr{H}},\widetilde{\mathbb E})$  as defined in Proposition \ref{lem3.4}. Denote  $S_n=\sum_{i=1}^n X_i$, $\tilde S_n=\sum_{i=1}^n \tilde X_i$.   Assume that each $X_n$ is tight. Consider the following statements:
\begin{description}
  \item[\rm (i)]  There exists a $\mathcal F$-measurable finite random variables $S$     such that $S_n\to S$ a.s. $\outCapc$, i.e.,
  \begin{equation}\label{eqth3.1}
  \outCapc\left(\left\{\omega: \lim_{n\to \infty} S_n(\omega)\ne S(\omega)\right\}\right)=0;
  \end{equation}
  \item[\rm (ii)]  There exists a $\mathcal F$-measurable finite random variables $S$     such that $S_n \to S$ in $\outCapc$, i.e.,
  \begin{equation}\label{eqth3.2}
  \outCapc\left(|S_n-S|\ge \epsilon\right)\to 0 \text{ as } n\to\infty \; \text{ for all } \epsilon>0;
  \end{equation}
  \item[\rm (i$^{\prime}$)]  For the copy $\{\tilde X_n;n\ge 1\}$, there exists a $\sigma(\widetilde{\mathscr{H}})$-measurable finite random variables $\tilde S$     such that $\tilde S_n\to \tilde S$ a.s. $\widetilde{\Capc}^{\ast}$, i.e.,
  \begin{equation}\label{eqth3.1ad}
  \widetilde{\Capc}^{\ast}\left(\left\{\omega: \lim_{n\to \infty} \tilde S_n(\omega)\ne \tilde S(\omega)\right\}\right)=0;
  \end{equation}
  \item[\rm (ii$^{\prime}$)]  For the copy $\{\tilde X_n;n\ge 1\}$, there exists a $\sigma(\widetilde{\mathscr{H}})$-measurable finite random variables $\tilde S$     such that $\tilde S_n\to \tilde S$  in  $\widetilde{\Capc}^{\ast}$, i.e.,
  \begin{equation}\label{eqth3.2}
  \widetilde{\Capc}^{\ast}\left(|\tilde S_n-\tilde S|\ge \epsilon\right)\to 0 \text{ as } n\to\infty \; \text{ for all } \epsilon>0;
  \end{equation}
  \item[\rm (iii)]  $\{S_n\}$ is a Cauchy sequence under $\upCapc$, i.e.,
  \begin{equation}\label{eqth3.3}
  \upCapc\left(|S_n-S_m|\ge \epsilon\right)\to 0 \text{ as } n,m\to\infty \; \text{ for all } \epsilon>0;
  \end{equation}
  \item[\rm (iv)] For some (equivalently, for any) $c>0$,
  \begin{description}
  \item[\rm (S1) ] $\sum\limits_{n=1}^{\infty}\upCapc(|X_n|>c)<\infty$,
  \item[\rm (S2) ] $\sum\limits_{n=1}^{\infty}\Sbep[X_n^{(c)}]$ and  $ \sum\limits_{n=1}^{\infty}\Sbep[-X_n^{(c)}]$ are both convergent,
  \item[\rm (S3) ]  $\sum\limits_{n=1}^{\infty}\Sbep\left[ \big(X_n^{(c)}-\Sbep[X_n^{(c)}]\big)^2\right] <\infty$ or/and $\sum\limits_{n=1}^{\infty}\Sbep\left[ \big(X_n^{(c)}+\Sbep[-X_n^{(c)}]\big)^2\right] <\infty$.
\end{description}
  \item[\rm (v)] $S_n$ converges in distribution, that is,  there is a sub-linear space $(\bar{\Omega}, \bar{\mathscr{H}}, \bar{\mathbb E})$ and a random variable $\bar{S}$ on it such that $\bar{S}$ is tight under $\bar{\mathbb E}$, i.e., $\bar{\Capc}(|\bar{S}|> x)\to 0$ as $x\to \infty$, and
       \begin{equation}\label{eqth3.4} \Sbep\left[\phi(S_n)\right]\to \bar{\mathbb E} \left[\phi(\bar{S})\right],\;\; \phi\in C_{b,Lip}(\mathbb{R}).
 \end{equation}
 \end{description}
 Then (i$^{\prime}$), (ii$^{\prime}$), (iii)-(v) are equivalent and each of them implies (i) and (ii). Further, suppose the space $(\Omega,\mathscr{H},\Sbep)$ satisfies one of the conditions (a)-(d) in Lemma \ref{lemBC2}. Then (i$^{\prime}$), (ii$^{\prime}$), (i)-(v) are equivalent.
\end{theorem}
\begin{remark} In the theorem $\upCapc$ can be replaced by any a capacity $\Capc$ with the property \eqref{eq1.4} by \eqref{eqV-V}, $\widetilde{\mathbb V}^{\ast}$ can be replaced $\widetilde{\mathbb C}^{\ast}$ or $\widetilde{V}$, $\outCapc$ can be replaced by $\mathbb C^{\ast}$, and, when one of the conditions (a)-(d) in Lemma \ref{lemBC2} is satisfied, $\outCapc$ can be replaced by $\Capc^{\mathscr{P}}$.
\end{remark}

At last, we give  an analogues  of Theorem  \ref{thLLN1}  for random vectors. Now, let $\{\bm X_n;n\ge 1\}$ be a sequence of i.i.d. random vectors in a sub-linear expectation space $(\Omega,\mathscr{H},\Sbep)$ who take values in an Euclidean space $\mathbb R^d$ with norm $|\bm x|=\sqrt{\sum_{i=1}^d x_i^2}$, and $\bm X_i\overset{d}=\bm X$. Suppose
\begin{equation} \label{eqthLLNv.0}
\lim_{c,d\to \infty}\Sbep\left[(|\bm X|\wedge d-c)^+\right]=0.
\end{equation}
 Then by Proposition \ref{prop1.1}, for any $\bm p\in \mathbb R^d$, $\vSbep[\langle \bm p,\bm X\rangle]$ is well-defined and finite, and
$$g(\bm p)=: \vSbep[\langle \bm p,\bm X\rangle], \;\; \bm p\in \mathbb R^d$$
is a sub-linear function defined on $\mathbb R^d$.  The assumption \eqref{eqthLLNv.0} is implied by a strong one as
\begin{equation}\label{eqChoquet1}
C_{\upCapc}(|\bm X|)<\infty
\end{equation}
or
$\lim_{c\to  \infty}\Sbep\left[(|\bm X|-c)^+\right]=0.
$
Further, if $\lim_{c\to  \infty}\Sbep\left[(|\bm X|-c)^+\right]=0$, then $g(\bm p)=  \Sbep[\langle \bm p,\bm X\rangle]$.  For the sub-linear function $g(\bm p)$, by Theorem 1.2.1 of Peng (2019),
   there exists     a (unique) bounded, convex and closed subset $\mathbb M$  such that (c.f., Peng (2019,Page 32))
$$g(\bm p)=\vSbep[\langle \bm p,\bm X\rangle]=   \sup_{\bm x\in\mathbb M}\langle \bm p,\bm x\rangle, \;\; \bm p\in \mathbb R^d.$$
We denote this set $\mathbb M$ by $\mathbb M_{\bm X}$ or $\mathbb M[\bm X]$. If $\bm X$ is a one-dimensional random variables, then $\mathbb M[\bm X]=[\cSbep[X], \Sbep[X]]$. For the multi-dimension case, recall $ \vSbep[\bm X]=(\vSbep[X_1],\ldots,\vSbep[X_d])$ and $E[\bm X] =(E[X_1],\ldots,E[X_d])$ for $\bm X=(X_1,\ldots,X_d)$.
\begin{lemma}   Under the condition \eqref{eqthLLNv.0} we have
\begin{equation} \label{eqthLLNM=M}  \mathbb M_{\bm X}= \widetilde{\mathbb M}_{\bm X}=:\Big\{ E[\bm X]: E\in \mathscr{E} \Big\}, \; \text{ where } \mathscr{E}  \text{ is defined as \eqref{linearexpression2.1}.}
 \end{equation}
 \end{lemma}
{\bf Proof}. It is obvious that
$$\sup_{\bm x\in \widetilde{\mathbb M}_{\bm X}}\langle \bm p,\bm x\rangle
 =\sup_{E\in\mathscr{E}}E[\langle \bm p,\bm X\rangle]=\vSbep[\langle \bm p,\bm X\rangle]=\sup_{\bm x\in  \mathbb M_{\bm X}}\langle \bm p,\bm x\rangle\;\; \text{ for all } \bm p\in \mathbb R^d.
 $$
 For \eqref{eqthLLNM=M}, it is sufficient to show that $\widetilde{\mathbb M}_{\bm X}$ is also a bounded, convex and closed subset of $\mathbb R^d$. The boundedness and convexity are obvious.   Next, we show that it is  closed. Suppose $E_i\in \mathscr{E}$, $E_i[\bm X]\to \bm b$. We want to show that $\bm b\in \widetilde{\mathbb M}_{\bm X}$. For each $E_i$, define $\widetilde{E}_i$ by $\widetilde{E}_i[\varphi(\bm x)]=E_i[\varphi(\bm X)]$, $\varphi\in C_{l,Lip}(\mathbb R^d)$. It is easily checked that
 if $C_{b,Lip}(\mathbb R^d)\ni \varphi_n\searrow 0$, then
\begin{align*}
0\le & \widetilde{E}_i[\varphi_n(\bm x)]=E_i[\varphi_n(\bm X)]\le  \vSbep[\varphi_n(\bm X)]\\
\le & \sup_{|\bm x|\le c}|\varphi_n(\bm x)|+\|\varphi_1\|c^{-1} \vSbep[|\bm X|]\to 0,
\end{align*}
as $n\to \infty$ and then $c\to \infty$.
 By Daniell-Stone's theorem,  there exists a probability measure $P_i$ on $\mathbb R^d$ such that
 $$ E_i[\varphi(\bm X)]=\widetilde{E}_i[\varphi(\bm x)]=P_i[\varphi(\bm x)],\; \; \forall \;  \varphi\in C_{b,Lip}(\mathbb R^d). $$
 Note that $\sup_iP_i(|\bm x|\ge c)\le c^{-1} \sup_iP_i[|\bm x|]\le c^{-1}\vSbep[|\bm X|]\to 0 \text{ as } c\to \infty$.
 So, on $\mathbb R^d$, the sequence $\{P_i\}$ is tight and so is relatively weakly compact. Then, there exist a subsequence $i_j$ and a probability $P$ on $\mathbb R^d$ such that
 \begin{equation}\label{eqConvofE} E_{i_j}[\varphi(\bm X)]=P_{i_j}[\varphi(\bm x)]\to P[\varphi(\bm x)]\; \;\forall \varphi\in C_{b,Lip}(\mathbb R^d).
 \end{equation}
 On the   the space $\mathscr{L}=\{Y=\varphi(\bm X): \varphi\in C_{l,Lip}(\mathbb R^d),Y\in \mathscr{H}_1\}$ we define an operator $E$ by
 $$ E[Y]=\lim_{j\to \infty} E_{i_j}[Y],\;\; Y\in \mathscr{L}. $$
 First, by \eqref{eqConvofE}, $E$ is well defined for bounded $Y\in \mathscr{L}$. Note
 $$ |E_{i_j}[Y]-E_{i_j}[(-c)\vee Y\wedge c]|=|E_{i_j}[Y-(-c)\vee Y\wedge c]|\le \vSbep[(|Y|-c)^+]\to 0 \text{ as } c\to\infty $$
 for $Y\in \mathscr{L}$. $E[Y]$ is well defined on $\mathscr{L}$ and $E[Y]=\lim_{c\to\infty}E[(-c)\vee Y\wedge c]$. It follows that
 $$b=\lim_{j\to \infty}E_{i_j}[\bm X]=E[\bm X]. $$
Since each $E_{i_j}\in\mathscr{E}$ is a finite additive linear expectation with $E_{i_j}\le \vSbep$, its limit $E$ is also a finite additive linear expectation on $\mathscr{L}$ with $E\le \vSbep$.  By the Hahn-Banach theorem, there exists  a finite additive linear expectation  $E^e$ defined on $\mathscr{H}_1$ such that, $E^e=E$ on $\mathscr{L}$ and, $E^e\le \vSbep$ on $\mathscr{H}_1$. So $E^e\in\mathscr{E}$. Hence, $b=E[\bm X]=E^e[\bm X]\in \widetilde{\mathbb M}_{\bm X}$. It follows that  $\widetilde{\mathbb M}_{\bm X}$ is a closed. \eqref{eqthLLNM=M} is proved. $\Box$

\bigskip

The following is the strong law of large numbers for i.i.d. random vectors.  Let $\{\bm X_n;n\ge 1\}$ be a sequence of i.i.d. random variables in a sub-linear expectation space $(\Omega,\mathscr{H},\Sbep)$, $\bm X_i\overset{d}= \bm X$.  Denote $\bm S_n=\sum_{i=1}^n \bm X_i$.

\begin{theorem}\label{thLLNv}
  If \eqref{eqChoquet1} is satisfied,
  then
  \begin{equation}\label{eqthLLNv.2}
  \outcCapc\left(C\left\{\frac{\bm S_n}{n}\right\}\subset \mathbb M_{\bm X} \right)=1.
  \end{equation}
  Further, suppose the space $(\Omega,\mathscr{H},\Sbep)$ satisfies one of the conditions (a)-(d) in Lemma \ref{lemBC2}. Then for $\Capc=\Capc^{\mathscr{P}}$, $\mathbb C^{\ast}$ or $\outCapc$,
  \begin{equation}\label{eqthLLNv.3}
 \Capc\left(C\left\{\frac{ \bm S_n}{n}\right\}= \mathbb M_{\bm X} \right)=1.
  \end{equation}
Further,
\begin{equation}\label{eqcorLLNv.4}
   \Capc\left(\lim_{n\to \infty} \frac{\bm S_n}{n}=\bm b \right)= \begin{cases} 1, & \text{ when } b\in \mathbb M_{\bm X}, \\
  0, & \text{ when } b\not\in \mathbb M_{\bm X}.
  \end{cases}
  \end{equation}
\end{theorem}

\eqref{eqthLLNv.3}   tells us that, under the upper capacity,  the limits of $\frac{\bm S_n}{n}$   fills the set $\mathbb M_{\bm X}$.
The following corollary  tells that, under lower capacity, the limit  of $\frac{ \bm S_n}{n}$ can only be a point.
\begin{corollary}\label{cor3.3} Suppose the space $(\Omega,\mathscr{H},\Sbep)$ satisfies one of the conditions (a)-(d) in Lemma \ref{lemBC2}.
 Assume that \eqref{eqChoquet1}  is satisfied.  If $\Capc=\Capc^{\mathscr{P}}$, $\mathbb C^{\ast}$ or $\outCapc$, there exists a subset $\mathbb O$ of $\mathbb R^d$ such that
  \begin{equation}\label{eqcor3.3.1}
\cCapc\left(C\left\{\frac{\bm S_n}{n}\right\}= \mathbb O \right)>0,
  \end{equation}
 then
  \begin{equation} \label{eqcor3.3.2}
  \vSbep[-\bm X]=-\vSbep[\bm X]\;  \text{ and } \; \mathbb O=\{\vSbep[\bm X]\}.
  \end{equation}
 \end{corollary}

 This is a direct corollary of Theorem \ref{thLLNv}. In fact, combining \eqref{eqcorLLNv.4}  and  \eqref{eqcor3.3.1} yields
 $$\Capc\left(\lim_{n\to \infty} \frac{ \bm S_n}{n}=\bm b \text{ and } C\left\{\frac{ \bm S_n}{n}\right\}= \mathbb O\right)>0  \text{ for all  } b\in \mathbb M_{\bm X}. $$
It follows that $\mathbb O=\{\bm b\}$ for all $b\in \mathbb M_{\bm X}$. Hence $\mathbb M_{\bm X}$ has only one point and then
\eqref{eqcor3.3.2} holds.

To prove Theorem \ref{thLLNv}, we need a weak law of large number which is of independent interest.
\begin{proposition}\label{lemWLLN}
 Let $\{\bm X_n;n\ge 1\}$ be a sequence of i.i.d. random variables in a sub-linear expectation space $(\Omega,\mathscr{H},\Sbep)$, $\bm S_n=\sum_{i=1}^n \bm X_i$.
  If $\lim\limits_{c,d\to\infty}\Sbep[(|\bm X_1|\wedge d-c)^+]=0$, then
  \begin{equation}\label{eqlemWLLN1}
  \upCapc\left(\frac{\bm S_n}{n}\not\in \mathbb M_{\bm X}^{\epsilon}\right)=\upCapc\left(dist\big( \bm S_n/{n}, \mathbb M_{\bm X} \big)\ge \epsilon\right)\to 0 \text{ for all } \epsilon>0
  \end{equation}
  and
   \begin{equation}\label{eqlemWLLN2}
  \upCapc\left(\Big|\frac{\bm S_n}{n}-\bm b\Big|<\epsilon\right)\to 1 \text{ for all } \bm b\in \mathbb M_{\bm X} \text{ and } \epsilon>0,
  \end{equation}
 where $dist(\bm y,\mathbb M_{\bm X})=\inf\{|\bm y-\bm x|:\bm x\in \mathbb M_{\bm X}\}$,
$\mathbb M_{\bm X}^{\epsilon}=\{\bm y: |\bm y-\bm x|<\epsilon \text{ for some } \bm x \in \mathbb M_{\bm X}\}$ is the $\epsilon$-neighborhood of $\mathbb M_{\bm X}$. In particular,
   \begin{equation}\label{eqlemWLLN3} \lim_{n\to \infty}\Sbep\left[\varphi\left(\frac{\bm S_n}{n}\right)\right]=\sup_{\bm x\in \mathbb M_{\bm X}}\varphi(\bm x), \;\; \text{ for all }\; \varphi\in C_{b,Lip}(\mathbb R^d).
    \end{equation}
\end{proposition}

The weak law of large numbers \eqref{eqlemWLLN3} is proved by Peng (2019) under the condition that $\Sbep[(|\bm X_1|-c)^+]\to 0$ as $c\to \infty$, by considering the solutions of the following parabolic PDEs defined on $[0,\infty)\times \mathbb R^d$,
$$ \partial_t u -g(D u)=0, \;\; u\big|_{t=0}=\varphi.$$
 \iffalse Now, we consider the sub-linear expectation $(\Omega,\mathscr{H}_1,\vSbep)$ instead, where $\mathscr{H}_1$ is defined as in Proposition \ref{prop1.1}.
Note $\vSbep[X]=\Sbep[X]$ if $X$ is bounded, and $\vSbep[(|X_1|-c)^+]\to 0$ as $c\to\infty$. \eqref{eqlemWLLN3}   can follow directly.
\fi
For the completeness of this paper, we will give a purely probabilistic proof in which only the probability inequalities are used.

 \section{Proofs of the law of large numbers}\label{sectProof}
\setcounter{equation}{0}
Before the proofs, we need a more lemma.
 \begin{lemma}\label{lem3.6} Suppose $X\in \mathscr{H}$, $1\le p<2$, $C_{\upCapc}(|X|^p)<\infty$. Then
 \begin{equation}\label{eqlem3.6.1}
 \sum_{i=1}^{\infty} \upCapc\left(|X|\ge M i^{1/p}\right)<\infty, \; \forall M>0,
 \end{equation}
  \begin{equation}\label{eqlem3.6.2}
 \sum_{i=1}^{\infty} \frac{\Sbep[X^2\wedge (Mi^{2/p})]}{i^{2/p}}<\infty, \; \forall M>0
 \end{equation}
 and
 \begin{equation}\label{eqlem3.6.3}
  \vSbep\left[(|X|-c)^+\right]=o(c^{1-p}) \text{ and } \Sbep[X^2\wedge c^2]=o(c^{2-p})\; \text{ as } c\to\infty.
 \end{equation}
 Further,
  \begin{equation}\label{eqlem3.6.4}
  C_{\upCapc}(|X|^p)=\infty\Longleftrightarrow \sum_{i=1}^{\infty} \upCapc\left(|X|\ge M i^{1/p}\right)=\infty, \; \forall M>0.
   \end{equation}
 \end{lemma}

 {\bf Proof}. \eqref{eqlem3.6.1} and \eqref{eqlem3.6.4}  are obvious by noting $C_{\upCapc}(|X|^p)=\int_0^{\infty} \upCapc\left(|X|>x^{1/p}\right)dx$. \eqref{eqlem3.6.2}  is similar to Lemma 3.9 (a) of Zhang (2016) and is proved in Zhang and Lin (2018). For  \eqref{eqlem3.6.3}, we have
 \begin{align*}
 &\vSbep\left[(|X|-c)^+\right]\le C_{\upCapc}\left((|X|-c)^+\right)=\int_c^{\infty}\upCapc(|X|>x)dx \\
 =&\frac{1}{p}\int_{c^p}^{\infty}y^{1/p-1}\upCapc(|X|^p>y)dy\le \frac{1}{p}c^{1-p}\int_{c^p}^{\infty} \upCapc(|X|^p>y)dy=o(c^{1-p})
 \end{align*}
 and
 \begin{align*}
 &\vSbep\left[X^2\wedge c^2\right]\le C_{\upCapc}\left(X^2\wedge c^2\right)=\int_0^{c^2}\upCapc(X^2>x)dx \\
 =&\frac{2}{p}\int_0^{c^p}y^{2/p-1}\upCapc(|X|^p>y)dy =o(c^{2-p}).
 \end{align*}
The proof is completed. $\Box$

\subsection{One-dimensional case}
 Now we turn to the proofs of the main results. We first consider the LLN for one-dimensional random variables.

 {\bf Proof of Theorems \ref{thLLN1} and \ref{thLLN2}}.  When \eqref{eqthLLN1.1} is satisfied, each $X_n$ is tight. Obviously, \eqref{eqthLLN1.4} is implied by \eqref{eqthLLN1.3} by noting that $\outcCapc\left(\frac{S_n}n-\frac{S_{n-1}}{n-1}\to 0\right)=1$.  \eqref{eqthLLN1.2} and \eqref{eqthLLN1.3} are special cases of \eqref{eqthLLN2.2} and \eqref{eqthLLN2.3}, respectively.  For \eqref{eqthLLN2.2}, we let $Z_{k,i}=(-2^{k/p})\vee X_i\wedge 2^{k/p}$, $i=1,\ldots, 2^k$. Then
 $$\left|\Sbep[Z_{k,i}]-\vSbep[X_i]\right|\le \vSbep\left[(|X_1|-2^{k/p})^+\right]=o\left(2^{k(1/p-1)}\right) $$
 by Lemma \ref{lem3.6}. For any $\epsilon>0$, by Lemma \ref{lem4.1} and \eqref{eqV-V} we have for $k$ large enough,
 \begin{align*}
 &\upCapc\left(\max_{2^{k-1}\le n\le 2^k}\frac{S_n-n\vSbep[X_1]}{n^{1/p}}\ge \epsilon\right)\\
 \le & \upCapc\left(\max_{2^{k-1}\le n\le 2^k}\sum_{i=1}^n (X_i-\vSbep[X_i])\ge \epsilon 2^{(k-1)/p}\right)\\
\le & \upCapc\left(\max_{ n\le 2^k}\sum_{i=1}^n (Z_{k,i}-\vSbep[Z_{k,i}])\ge \epsilon 2^{k/p}/4\right)+\upCapc\left(\max_{i\le 2^k}|X_i|>2^{k/p}\right)\\
\le & C 2^{-2k/p}\sum_{i=1}^{2^k} \Sbep[Z_{k,i}^2]+\sum_{i=1}^{2^k}\upCapc\left(|X_1|>2^{k/p}/2\right)\\
=&C 2^{-2k/p}2^k \Sbep[X_1^2\wedge 2^{2k/p}]+2^k\upCapc\left(|X_1|>2^{k/p}/2\right)\\
\le &4C \sum_{i=2^{k}+1}^{2^{k+1}}\frac{\Sbep[X_1^2\wedge i^{2/p}]}{i^{2/p}}+2\sum_{i=2^{k-1}+1}^{2^k}\upCapc(|X_1|>i^{1/p}/2).
 \end{align*}
It follows that
\begin{align*}
 &\sum_{k=1}^{\infty}\outCapc\left(\max_{2^{k-1}\le n\le 2^k}\frac{S_n-n\vSbep[X_1]}{n^{1/p}}\ge \epsilon\right)\\
 \le &\sum_{k=1}^{\infty}\upCapc\left(\max_{2^{k-1}\le n\le 2^k}\frac{S_n-n\vSbep[X_1]}{n^{1/p}}\ge \epsilon\right)\\
\le &4C \sum_{i=1}^{\infty}\frac{\Sbep[X_1^2\wedge i^{2/p}]}{i^{2/p}}+2\sum_{i=1}^{\infty}\upCapc(|X_1|>i^{1/p}/2)<\infty,
 \end{align*}
 by Lemma \ref{lem3.6}.  By noting that $\outCapc$ is a countably sub-additive capacity and the Borel-Cantelli (Lemma \ref{lemBCdirect}), we have
$$\outCapc\left(\limsup_{n\to\infty}\frac{S_n-n\vSbep[X_1]}{n^{1/p}}\ge \epsilon\right)\le \outCapc\left(\max_{2^{k-1}\le n\le 2^k}\frac{S_n-n\vSbep[X_1]}{n^{1/p}}\ge \epsilon\; i.o.\right)=0. $$
By the countable sub-additivity of $\outCapc$ again,
$$ \outCapc\left(\limsup_{n\to\infty}\frac{S_n-n\vSbep[X_1]}{n^{1/p}} >0\right)
=\outCapc\left(\bigcup_{l=1}^{\infty}\left\{\limsup_{n\to\infty}\frac{S_n-n\vSbep[X_1]}{n^{1/p}}\ge \frac{1}{l}\right\}\right)=0. $$
For $-X_i$s, we have a similar result.  \eqref{eqthLLN2.2} is proved.

For \eqref{eqthLLN2.3}, it is sufficient to show that
\begin{equation}\label{eqproofthLLN1}
\Capc^{\mathscr{P}}\left(\liminf_{n\to\infty}\frac{\tilde S_n-n\vcSbep[X_1]}{n^{1/p}}\le 0 \;\text{ and }\; \limsup_{n\to\infty}\frac{\tilde S_n-n\vSbep[X_1]}{n^{1/p}}\ge 0\right)=1.
\end{equation}
Let $Y_{ni}=(-n^{1/p})\vee X_i\wedge n^{1/p}$, $i=1,\ldots,n$. Then $\widetilde{\mathbb M}[Y_{ni}]=[\cSbep[Y_{ni}],\Sbep[Y_{ni}]]$. By Lemmas \ref{moment_v} and \ref{lem3.6},
 \begin{align*}
 &\lowCapc\left(\sum_{i=1}^n (-Y_{ni}+\Sbep[Y_{ni}])\ge \epsilon n^{1/p}\right)\\
 =& \lowCapc\left(\sum_{i=1}^n (-Y_{ni}-\cSbep[-Y_{ni}])\ge \epsilon n^{1/p}\right)\le 2\frac{n\Sbep[X_1^2\wedge n^{2/p}]}{\epsilon^2 n^{2/p}}\to 0.
 \end{align*}
On the other hand, $n|\Sbep[Y_{ni}-\vSbep[X_1]|\le n \vSbep\left[(|X_1|-n^{1/p})^+\right]=o(n^{1/p})$ and
$$ \upCapc\left(Y_{ni}\ne X_i, \; i=1,\ldots,n\right)\le n\upCapc(|X_1|>n^{1/p})\to 0, $$
by Lemma \ref{lem3.6}. It follows that
$$ \lowCapc\left(\sum_{i=1}^n (-X_i+\Sbep[X_1])\ge 2\epsilon n^{1/p}\right) \to 0. $$
That is
$$\upCapc\left(\frac{S_n-n\vSbep[X_1]}{n^{1/p}}\ge -\epsilon\right)\to 1\; \text{ for all } \epsilon>0. $$
By considering $-X_i$s, similarly we have
$$\upCapc\left(\frac{-S_n+n\vcSbep[X_1]}{n^{1/p}}\ge -\epsilon\right)\to 1\; \text{ for all } \epsilon>0. $$
For $\epsilon_k=1/2^k$, $k=1,2,\ldots$, we can choose $n_k$ successively such that  $n_k\nearrow \infty$, $n_{k-1}/n_k^{1/p}\to 0$, and
$$ \upCapc\left(\frac{S_{n_k}-S_{n_{k-1}}-(n_k-n_{k-1})\vSbep[X_1]}{(n_k-n_{k-1})^{1/p}}\ge -\epsilon_k \right)\ge 1-\epsilon_k,$$
$$ \upCapc\left(-\frac{S_{n_k}-S_{n_{k-1}}-(n_k-n_{k-1})\vcSbep[X_1]}{(n_k-n_{k-1})^{1/p}}\ge -\epsilon_k \right)\ge 1-\epsilon_k.$$
It follows that
$$ \sum_{k=1}^{\infty}\upCapc\left(\frac{S_{n_k}-S_{n_{k-1}}-(n_k-n_{k-1})\vSbep[X_1]}{(n_k-n_{k-1})^{1/p}}\ge -\epsilon_k \right)=\infty,$$
$$ \sum_{k=1}^{\infty}\upCapc\left(-\frac{S_{n_k}-S_{n_{k-1}}-(n_k-n_{k-1})\vcSbep[X_1]}{(n_k-n_{k-1})^{1/p}}\ge -\epsilon_k \right)=\infty.$$
Let
$$ A=\left\{\frac{ S_{n_k}-  S_{n_{k-1}}-(n_k-n_{k-1})\vSbep[X_1]}{(n_k-n_{k-1})^{1/p}}\ge -\epsilon_k \; i.o.\right\},$$
$$ B=\left\{-\frac{  S_{n_k}- S_{n_{k-1}}-(n_k-n_{k-1})\vcSbep[X_1]}{(n_k-n_{k-1})^{1/p}}\ge -\epsilon_k \; i.o.\right\}.$$
  By the Borel-Cantelli lemma (Lemma \ref{lemBC2}), $\Capc^{\mathscr{P}}(AB)=1$.
On $AB$ and $C=\left\{\limsup_{n\to\infty}\frac{| S_n|}{n}<\infty\right\}$,
\begin{align*}
\limsup_{n\to \infty}\frac{  S_n-n\vSbep[X_1]}{n^{1/p}}\ge & \limsup_{k\to\infty}\frac{  S_{n_k}-  S_{n_{k-1}}-(n_k-n_{k-1})\vSbep[X_1]}{n_k^{1/p}}\\
\ge & \limsup_{k\to\infty}\frac{  S_{n_k}-  S_{n_{k-1}}-(n_k-n_{k-1})\vSbep[X_1]}{(n_k-n_{k-1})^{1/p}}\ge 0,\\
\limsup_{n\to \infty}\frac{-  S_n+n\vcSbep[X_1]}{n^{1/p}}\ge &
\limsup_{k\to\infty}\left(- \frac{  S_{n_k}-  S_{n_{k-1}}-(n_k-n_{k-1})\vcSbep[X_1]}{(n_k-n_{k-1})^{1/p}}\right)\ge 0.
\end{align*}
Note $\Capc^{\mathscr{P}}(ABC)\ge   \Capc^{\mathscr{P}}(AB)-\Capc^{\mathscr{P}}(C^c)=1-0=1$ by \eqref{eqthLLN1.2}. The proof of \eqref{eqproofthLLN1} is completed.

For Theorem \ref{thLLN1} (b), suppose $C_{\upCapc}(|X_1|)=\infty$. Then
$$ \sum_{n=1}^{\infty}\upCapc(|X_n|\ge Mn)\ge  \sum_{n=1}^{\infty}\upCapc(|X_1|\ge Mn/2)=\infty, \; \text{ for all } M>0. $$
So, there exists a sequence $1<M_n\nearrow \infty$ such that
$$ \sum_{n=1}^{\infty}\upCapc(|X_n|\ge M_n n) =\infty. $$
By the Borel-Cantelli Lemma (Lemma \ref{lemBC2}),
$$ \Capc^{\mathscr{P}}\left(| X_n|\ge M_n n  \; i.o.\right)=1. $$
On the event $\{| X_n|\ge M_n n  \; i.o.\}$, we have
$$ \infty=\limsup_{n\to \infty}\frac{|  X_n|}{n}\le 2\limsup_{n\to \infty}\frac{|  S_n|}{n}. $$
It follows that
\begin{equation}\label{eqproofthLLN2}
\Capc^{\mathscr{P}}\left(\limsup_{n\to \infty}\frac{|  S_n|}{n}=\infty\right)=1,
\end{equation}
which contradicts with \eqref{eqthLLN1.5}. The proof is now completed. $\Box$

\bigskip

 {\bf Proof of Corollaries  \ref{cor1} and \ref{cor2}}. It is sufficient to show Corollary \ref{cor2}.

 (a)    When $b\not\in[\vcSbep[X_1],\vSbep[X_1]]$, the conclusion \eqref{eqcor2.1} is obvious by \eqref{eqthLLN2.2}. If $b\in[\vcSbep[X_1],\vSbep[X_1]]$, then there exists an $\alpha\in [0,1]$ such that $b=\alpha\vSbep[X_1]+(1-\alpha)\vcSbep[X_1]$. Let $Y_i=(-i^{1/p})\vee X_i\wedge i^{1/p}$ and $\mu_{\alpha,i}= \alpha\Sbep[Y_i]+(1-\alpha)\cSbep[Y_i]$. Then $\sum_{i=1}^n |\mu_{\alpha,i}-b|\le \sum_{i=1}^n \vSbep\left[(|X_1|-i^{1/p})^+\right]=o(n^{1/p})$ and $\outCapc( X_i\ne   Y_i \; i.o.)=0$. So, it is sufficient to show that for each $\alpha\in [0,1]$,
 \begin{equation}\label{eqproofcor1}
\Capc^{\mathscr{P}}\left(\lim_{n\to\infty}\frac{\sum_{i=1}^n(  Y_i-\mu_{\alpha,i})}{n^{1/p}}=0\right)=1.
 \end{equation}
 For each $i$, by the expression \eqref{linearexpression}, there exist $\theta_{i,1},\theta_{i,2}\in \Theta$ such that
 $$ E_{\theta_{i,1}}[Y_i]=\Sbep[Y_i]\; \text{ and }\;   E_{\theta_{i,2}}[Y_i]=\cSbep[Y_i]. $$
 Define the linear operator $E_i=\alpha E_{\theta_{i,1}}+(1-\alpha)E_{\theta_{i,2}}$. Then
 $$ E_i[Y_i]=\mu_{\alpha,i} \; \text{ and }\; E_i\le \Sbep. $$
 Note that each $Y_n$ is tight.  By Proposition \ref{lem3.4}, there exist a copy $\{\tilde Y_n;n\ge 1\}$ on $(\widetilde{\Omega},\widetilde{\mathscr{H}},\widetilde{\mathbb E})$ of $\{Y_n;n\ge 1\}$  and a probability measure $Q$ on $\sigma(\tilde Y_1,\tilde Y_2,\ldots)$ such that
 such that $\{\tilde Y_n;n\ge 1\}$ is a sequence of independent random variables under $Q$,
$$
Q\left[\varphi(\tilde Y_i)\right]=E_i\left[\varphi(Y_i)\right]\; \text{ for all } \varphi\in C_{b,Lip}(\mathbb R),
$$
$$
Q\left[\varphi(\tilde Y_1,\ldots,\tilde Y_d)\right]\le \Sbep\left[\varphi(Y_1,\ldots,Y_d)\right]\; \text{ for all } \varphi\in C_{b,Lip}(\mathbb R^d)
$$
and
\begin{equation}\label{eqproofcor2}
\widetilde{v}(B)\le P(B)\le \widetilde{V}(B) \; \text{ for all } B\in\sigma(\tilde Y_1,\tilde Y_2,\ldots).
\end{equation}
Note $|E_i[Y_i^{(c)}]-E_i[Y_i]|\le \Sbep[(|Y_i|-c)^+]\to 0$ as $c\to \infty$. We have
\begin{align}\label{eqExpQ}
Q[Y_i]=& \lim_{c\to \infty}Q[Y_i^{(c)}]=\lim_{c\to \infty}E_i[Y_i^{(c)}]
=E_i[Y_i]=\mu_{\alpha,i}, \\
\label{eqVarQ}
Q[Y_i^2]=& \lim_{c\to \infty}Q[\tilde Y_i^2 \wedge c]=\lim_{c\to \infty} E_i[ Y_i^2\wedge c]\le \Sbep[Y_i^2].
\end{align}
Then
$$ \sum_{i=1}^{\infty}\frac{Q[\tilde Y_i^2]}{i^{2/p}}
\le \sum_{i=1}^{\infty}\frac{\Sbep[Y_i^2]}{i^{2/p}}=\sum_{i=1}^{\infty}\frac{\Sbep[X_1^2\wedge i^{2/p}]}{i^{2/p}}<\infty,
 $$
 by Lemma \ref{lem3.6}. Denote $  T_n=\sum_{i=1}^n(Y_i-\mu_{\alpha,i})$ and  $\tilde T_n=\sum_{i=1}^n(\tilde Y_i-\mu_{\alpha,i})$, $n_k=2^k$. Then
\begin{align*}
 &Q\left(\frac{\max_{n_k+1\le n\le n_{k+1}}|\tilde T_n-\tilde T_{n_k}|}{n_k^{1/p}}>   \epsilon  \right)\\
 \le &2\epsilon^{-2}n_k^{-2/p}\sum_{i=n_k+1}^{n_{k+1}}Q[\tilde Y_i^2] \le \epsilon^{-2} 2^{2/p+1}\sum_{i=n_k+1}^{n_{k+1}}\frac{Q[\tilde Y_i^2]}{i^{1/p}}.
 \end{align*}
  It follows that
  $$\sum_{k=1}^{\infty}Q\left(\frac{\max_{n_k+1\le n\le n_{k+1}}|\tilde T_n-\tilde T_{n_k}|}{n_k^{1/p}}>   \epsilon  \right)<\infty, \;\; \text{ for all } \epsilon>0. $$
  Then there exists a sequence $\epsilon_k\searrow 0$ such that
  $$\sum_{k=1}^{\infty}Q\left(\frac{\max_{n_k+1\le n\le n_{k+1}}|\tilde T_n-\tilde T_{n_k}|}{n_k^{1/p}}>   \epsilon_k  \right)<\infty.$$
  By \eqref{eqproofcor2} and \eqref{eqV-V},
  \begin{align}\label{eqconvergenceforBC}
  &\sum_{k=1}^{\infty}\cCapc^{\mathscr{P}}\left(\frac{\max_{n_k+1\le n\le n_{k+1}}|  T_n-  T_{n_k}|}{n_k^{1/p}}>   2\epsilon_k  \right)
  \nonumber \\
 \le  &\sum_{k=1}^{\infty}\widetilde{v}\left(\frac{\max_{n_k+1\le n\le n_{k+1}}|\tilde T_n-\tilde T_{n_k}|}{n_k^{1/p}}>   \epsilon_k  \right)<\infty.
  \end{align}
  Note the independence. By Lemma \ref{lemBC2},
  $$\cCapc^{\mathscr{P}}\left(A_k\; i.o.  \right)=0 \; \text{ with } \; A_k=\left\{\frac{\max_{n_k+1\le n\le n_{k+1}}|  T_n-  T_{n_k}|}{n_k^{1/p}}>   2\epsilon_k\right\}. $$
  Note that on the event $(A_k\; i.o.)^c$,
  $$ \lim_{k\to \infty}\frac{\max_{n_k+1\le n\le n_{k+1}}|  T_n-  T_{n_k}|}{n_k^{1/p}}=0, $$
  which implies $ \lim_{n\to \infty} \frac{T_n}{n}=0$. \eqref{eqproofcor1} is proved.

(b)  First, note the facts that  $\Capc(AB)=1$ whenever $\Capc(A)=\cCapc(B)=1$, $\cCapc(AB)=1$ whenever $\cCapc(A)=\cCapc(B)=1$.  If \eqref{eqthLLN2.1} and \eqref{eqcor2.3} hold,  and $\vcSbep[X_1]=\vSbep[X_1]$, then \eqref{eqcor2.2} is obvious by \eqref{eqthLLN2.2}.
Conversely, suppose  \eqref{eqcor2.2} holds. Let $A=\left\{\liminf\limits_{n\to \infty} \frac{  S_n-n\vcSbep[X_1]}{n^{1/p}}=0\right\}$, $B=\left\{\limsup\limits_{n\to \infty} \frac{  S_n-n\vSbep[X_1]}{n^{1/p}}=0\right\}$ and $C=\left\{\lim\limits_{n\to \infty} \frac{  S_n-nb}{n^{1/p}}=c\right\}$.

We first consider the case $p=1$. If $C_{\upCapc}(|X_1|)=\infty$, then \eqref{eqproofthLLN2} holds, which contradicts with \eqref{eqcor2.2}. So, $C_{\upCapc}(|X_1|)<\infty$. By \eqref{eqcor2.2}) and \eqref{eqthLLN2.3} with $p=1$, $\Capc(ABC)=1$.  While, on $ABC$,
\begin{align*}
\vcSbep[X_1]=& \liminf_{n\to \infty}\frac{  S_n}{n}=c+b, \\
\vSbep[X_1]=& \limsup_{n\to \infty}\frac{  S_n}{n}=c+b.
\end{align*}
It follows that $\vSbep[X_1]=\vcSbep[X_1]$.  Then, by the direct part,
$$ \cCapc\left(\lim_{n\to \infty}\frac{  S_n}{n}=\vSbep[X_1]\right)=1, $$
which, together with \eqref{eqcor2.2}, implies $\cCapc(b+c=\vSbep[X_1])=1$.

Now, suppose $1<p<2$. Then
$$\cCapc\left(\lim_{n\to \infty}\frac{ S_n}{n}=b\right)=1$$
by \eqref{eqcor2.2}. By the conclusion for the case $p=1$, we must have $\vSbep[X_1]=\vcSbep[X_1]$ and $\cCapc(b=\vcSbep[X_1])=1$. Suppose $C_{\upCapc}(|X_1|^p)=\infty$. Then $C_{\upCapc}(|X_1-\vSbep[X_1]|^p)=\infty$. Similar to \eqref{eqproofthLLN2}, we have
\begin{align*}
 & \Capc \left(\limsup_{n\to \infty}\frac{| S_n-n\vSbep[X_1]|}{n^{1/p}}=\infty\right)\ge \Capc^{\mathscr{P}}\left(\limsup_{n\to \infty}\frac{| S_n-n\vSbep[X_1]|}{n^{1/p}}=\infty\right)\\
 & \;\; \ge
 \Capc^{\mathscr{P}}\left(\limsup_{n\to \infty}\frac{|  X_n-\vSbep[X_1]|}{n^{1/p}}=\infty\right)=1,
\end{align*}
 which contradicts with \eqref{eqcor2.2} by noting $\cCapc(b=\vSbep[X_1])=1$. So, \eqref{eqthLLN2.1} holds.
 By the direct part, \eqref{eqcor2.1} holds. Then
$$ \cCapc\left(\lim_{n\to \infty}\frac{   S_n-nb }{n^{1/p}}=0\right)=
\cCapc\left(\lim_{n\to \infty}\frac{   S_n-n\vSbep[X_1] }{n^{1/p}}=0\right)=1, $$
 which, together with \eqref{eqcor2.2}, implies $\cCapc(c=0)=1$. The proof is now completed. $\Box$

\bigskip

Next, we consider the convergence of infinite series.

 {\bf Proof of Theorem \ref{th3}}.   (iii)$\Longleftrightarrow$(iv) and (v)$\implies$(iii) are proved by Zhang (2020) (c.f. Theorem 3.2, Theorem 3.3 there).

 First, we show that (iii)$\implies$(v). Let $\bar{\Omega}=\mathbb R$, $\bar{\mathscr{H}}=C_{b,Lip}(\mathbb R)$. Define $\bar{\mathbb E}$ by
 $$ \bar{\mathbb E}[\varphi]=\limsup_{n\to \infty} \Sbep[\varphi(S_n)], \;\; \varphi\in \bar{\mathscr{H}}, $$
 and define the random variable $\bar{S}$ by $\bar{S}(x)=x$.
 Since  each $X_i$ is tight,   each $S_n$ is tight. By \eqref{eqth3.3}, we have $\lim_{c\to \infty}\max_n\upCapc(|S_n|>c)=0$.  Then
 $$ \bar{\Capc}(|\bar{S}|\ge c)\le \limsup_{n\to \infty}\upCapc(|S_n|\ge c/2)\to 0 \text{ as } c\to \infty. $$
It follows that $\bar{S}$ is tight. For $\varphi\in \bar{\mathscr{H}}$, let $\varphi_c(x)=\varphi((-c)\vee x\wedge c)$. Then $\varphi_c$ is a uniformly continuous function. For any $\epsilon>0$, there is a $\delta>0$ such that $|\varphi_c(x)-\varphi_c(y)|<\epsilon$ when $|x-y|<\delta$. Hence
$$ \big|\Sbep[\varphi_c(S_n)]-\Sbep[\varphi_c(S_m)]\big|<\epsilon+\|\varphi|\upCapc(|S_n-S_m|>\delta)\to \epsilon $$
as $n,m\to \infty$. Hence $\Sbep[\varphi_c(S_n)]$ converges. It follows that
$$ \bar{E}[\varphi_c(\bar{S})]=\lim_{n\to \infty} \Sbep[\varphi_c(S_n)]. $$
On the other hand,
$$\max_n\big|\Sbep[\varphi_c(S_n)]-\Sbep[\varphi(S_n)]\big|\le \|\varphi\|\max_n \upCapc(|S_n|>c)\to 0 $$
and
$$\max_n\big|\bar{\mathbb E}[\varphi_c(\bar S)]-\bar{\mathbb E}[\varphi(\bar S)]\big|\le \|\varphi\|  \bar{\Capc}(|\bar S |>c)\to 0, $$
as $c\to \infty$. It follows that
$$ \bar{\mathbb E}[\varphi(\bar S)]=\lim_{n\to \infty} \Sbep[\varphi(S_n)], ;\; \forall \varphi\in C_{b,Lip}(\mathbb R). $$
(v) holds.

 Next, we show (iii)$\implies$ (i) and (ii).   By the L\'evy inequality \eqref{eqLIQ2}, it follows from  \eqref{eqth3.3} that
\begin{equation}\label{eqproofth3.1}
\upCapc\left(\max_{m\le i\le n}|S_i-S_m|\ge \epsilon\right)\to 0 \text{ as } n,m\to\infty \; \text{ for all } \epsilon>0.
\end{equation}
Let $\epsilon_k=1/2^k$. There exists a sequence $n_k\nearrow \infty$ such that
$$ \outCapc\left(\max_{n_k\le i\le n_{k+1}}|S_i-S_{n_k}|\ge \epsilon_k\right)\le \upCapc\left(\max_{n_k\le i\le n_{k+1}}|S_i-S_{n_k}|\ge \epsilon_k\right)<\epsilon_k. $$
It follows that
$$ \sum_{k=1}^{\infty}\outCapc\left(\max_{n_k\le i\le n_{k+1}}|S_i-S_{n_k}|\ge \epsilon_k\right)<\sum_{k=1}^{\infty}\epsilon_k<\infty.  $$
Note the countable sub-additivity of $\outCapc$. By the Borel-Cantelli lemma (Lemma \ref{lemBCdirect}),
$$ \outCapc(A)=0 \; \text{ where } A=\left\{\max_{n_k\le i\le n_{k+1}}|S_i-S_{n_k}|\ge \epsilon_k\; i.o.\right\}.$$
on $A^c$, $S=S_{n_0}+\sum_{k=1}^{\infty}(S_{n_k}-S_{n_{k-1}})$ is finite. Let $S(\omega)=0$ when $\omega\in A$.
On $A^c$, $S_{n_k}\to S$  and $\max_{n_k\le i\le n_{k+1}}|S_i-S_{n_k}|\to 0$ as $k\to\infty$, and so $S_i\to S$ as $i\to\infty$. Then (i) is proved.

Also,   on the event $\bigcap_{m=k}^{\infty}\left\{\max_{n_m\le i\le n_{m+1}}|S_i-S_{n_m}|\le \epsilon_m\right\}$,
$$ |S-S_{n_k}|\le \sum_{m=k}^{\infty} |S_{n_{m+1}}-S_{n_m}|\le \sum_{m=k}^{\infty} 2^{-m}=2^{-k+1}.$$
It follows that
$$\outCapc\left(|S_{n_k}-S|>2^{-k+1}\right)\le \sum_{m=k}^{\infty}\upCapc\left(|S_{n_{m+1}}-S_{n_m}|>\epsilon_m\right)\le \sum_{m=k}^{\infty}\epsilon_m<2^{-k+1}. $$
On the other hand, for any $\epsilon>0$,  when $k$ is large enough such that $2^{-k+1}<\epsilon/2$,
$$\outCapc\left(|S_n-S|>\epsilon/2\right)\le \outCapc\left(|S_{n_k}-S|>2^{-k+1}\right)+\outCapc\left(|S_{n_k}-S_n|>\epsilon/2\right)\to 0, $$
as $n,k\to\infty$.  Then (ii) is proved.

Note \eqref{eqV-V} and $(X_1,\ldots,X_n)\overset{d}=(\tilde X_1,\ldots, \tilde X_n)$, $n\ge 1$. (iii) is equivalent to that it holds for $\tilde S_n$. So, it implies (i$^{\prime}$) and (ii$^{\prime}$).

Note the copy space $(\widetilde{\Omega},\widetilde{\mathscr{H}},\widetilde{\mathbb E})$ satisfies the condition (d) in Lemma \ref{lemBC2}. At last, it is sufficient to show  (i)$\implies$ (iii), and (ii) $\implies$ (iii),  when one of the conditions (a)-(b) in Lemma \ref{lemBC2} is satisfied.

   Suppose (iii) does not hold. Then there exist constants $\epsilon_0>0$, $\delta_0>0$ and sequence $\{m_k\}$ and $\{n_k\}$ with $m_k<n_k\le m_{k+1}$ such that
$$ \upCapc\left(|S_{n_k}-S_{m_k}|\ge \epsilon_0\right)\ge \delta_0. $$
Note the independence of $\{S_{n_k}-S_{m_k}; k\ge 1\}$ and
$$ \sum_{k=1}^{\infty} \upCapc\left(|S_{n_k}-S_{m_k}|\ge \epsilon_0\right)=\infty. $$
By the Borel-Cantelli lemma (Lemma \ref{lemBC2}),
$$\Capc^{\mathscr{P}}\left(\limsup_{k\to\infty} |  S_{n_k}-  S_{m_k}|\ge \epsilon_0/2\right)=1. $$
However, on the event $\{\lim_{n\to \infty}  S_n=  S\}$ we have $\limsup_{k\to \infty}|  S_{n_k}- S_{m_k}|=0$. Thus,
$$\Capc^{\mathscr{P}}\left(\left\{\omega:\lim _{n\to \infty}  S_n(\omega)\ne  S(\omega)\right\}\right)=1, $$
which contradicts with \eqref{eqth3.1ad}. So, (i)$\implies$(iii) is proved.

Now, suppose $\Capc^{\mathscr{P}}(|S_n-S|>\epsilon)\to 0$ for all $\epsilon>0$. Then
$$\Capc^{\mathscr{P}}(|S_n-S_m|>\epsilon)\to 0 \text{ as } n,m \to \infty, \forall \epsilon>0, $$
which is equivalent to \eqref{eqth3.3} by \eqref{eqV-V}, since both $\Capc^{\mathscr{P}}$ and $\upCapc$ have the property \eqref{eq1.4} and $S_n-S_m\in \mathscr{H}$.      The proof is completed.
$\Box$

 \subsection{Multi-dimensional case}

 Now, we consider the LLN for random vectors.

{\bf Proof of Proposition \ref{lemWLLN}}. Recall  $\bm X^{(c)}=(X_1^{(c)},\ldots,X_d^{(c)})$ and $X_i^{(c)}=(-c)\vee X_i\wedge c$ for $\bm X=(X_1,\ldots, X_d)$. Note
 \begin{align*}
  \upCapc\left(\left| \frac{\bm S_n}{n} -  \frac{\sum_{i=1}^n\bm X_i^{(c)}}{n} \right|\ge \epsilon\right)
 \le  \epsilon^{-1}  \vSbep\big[|\bm X_1-\bm X_1^{(c)}|\big]\to 0 \text{ as } c\to \infty
 \end{align*}
 and
\begin{align*}
\sup_{E\in \mathscr{E}} \Big| E[\bm X] - E[\bm X^{(c)}] \Big|\le  \vSbep\big[|\bm X_1-\bm X_1^{(c)}|\big]\to 0 \text{ as } c\to \infty.
  \end{align*}
Hence, without loss of generality we can assume $|\bm X_i|\le c$ and $|\bm X|\le c$.
 Let $\delta=\epsilon^2/(4c)$, and $\mathcal{N}_{\delta}=\{\bm p_1,\ldots,\bm p_K\}\subset \{\bm p: |\bm p|\le 2c\}$ be a $\delta$-net of $\{\bm p: |\bm p|\le 2c\}$.
 We have the following fact,
\begin{equation}\label{eqprooflemWLLN.1}  \bm y\not\in \mathbb M_{\bm X}^{\epsilon} \text{ and } |\bm y|\le c \Longrightarrow \langle \bm p_i,\bm y\rangle\ge \vSbep[
 \langle \bm p_i,\bm X\rangle]+\epsilon^2/2 \text{ for some } \bm p_i\in \mathcal{N}_{\delta}.
 \end{equation}
 In fact, for $\bm y\not\in \mathbb M_{\bm X}^{\epsilon}$, there exists $\bm o=E[\bm X]\in \mathbb M_{\bm X}$ such that $\tau=:\inf_{x\in\mathbb M_{\bm X}}|\bm y-\bm x|=|\bm y-\bm o|\ge \epsilon$. Let $\bm p=\bm y-\bm o$. Then $|\bm p|\le 2c$ and
 $\langle \bm p,\bm y\rangle=\langle \bm p,\bm o\rangle +\tau^2$. For any $\bm x\in\mathbb M_{\bm X}$ and $0\le \alpha\le 1$, $\bm z=\alpha \bm x+(1-\alpha)\bm o\in \mathbb M_{\bm X}$.
 Then
 $$ |\bm y-\bm o|^2\le |\bm y-\bm z|^2 =|\bm y-\bm o|^2+\alpha^2|\bm x-\bm o|^2+2\alpha \langle \bm x-\bm o,\bm o-\bm y\rangle\;\text{ for all } \alpha\in [0,1]. $$
 It follows that
 $$ \langle \bm p,\bm o\rangle-\langle \bm p,\bm x\rangle=\langle \bm x-\bm o,\bm o-\bm y\rangle\ge 0.
 $$
 So $\langle \bm p,\bm o\rangle\ge \langle \bm p,\bm x\rangle$. It follows that $\langle \bm p,\bm o\rangle\ge \sup_{\bm x\in \mathbb M_{\bm X}}\langle \bm p,\bm x\rangle=\vSbep[
 \langle \bm p,\bm X\rangle]$. It follows that $\langle \bm p,\bm y\rangle\ge \vSbep[
 \langle \bm p,\bm X\rangle]+\epsilon^2$.  Further,   for the $\bm p$, there exists a $\bm p_i\in \mathcal{N}_{\delta}$ such that $|\bm p-\bm p_i|<\delta$. Then
 $$ \langle \bm p_i,\bm y\rangle- \vSbep[
 \langle \bm p_i,\bm X\rangle]\ge \langle \bm p,\bm y\rangle- \vSbep[
 \langle \bm p,\bm X\rangle]-|\bm p_i-\bm p||\bm y| - |\bm p_i-\bm p|\vSbep[|\bm X|]\ge \epsilon^2/2. $$
  Hence \eqref{eqprooflemWLLN.1} follows.  Now,
 it follows from the inequality \eqref{eqlem4.1.2} that
 \begin{align*}
   \upCapc\left(\frac{\bm S_n}{n}\not\in\mathbb M_{\bm X}^{\epsilon}\right)
 \le &  \sum_{i\in \mathcal{N}_{\delta}} \upCapc\left(\langle \bm p_i,\bm S_n/n\rangle\ge \vSbep[
 \langle \bm p_i,\bm X\rangle]+\epsilon^2/2\right)\\
 = &   \sum_{i\in \mathcal{N}_{\delta}} \upCapc\left(\sum_{k=1}^n\big(\langle \bm p_i,\bm X_k\rangle- \Sbep[
 \langle \bm p_i,\bm X_k\rangle]\big)\ge n\epsilon^2/2\right)\\
 \le &2(e+1)\sum_{i\in \mathcal{N}_{\delta}}\frac{n\Sbep[ \langle \bm p_i,\bm X\rangle^2]}{\epsilon^4n^2/4}\to 0.
 \end{align*}
 The proof of \eqref{eqlemWLLN1} is completed.

For \eqref{eqlemWLLN2}, we suppose $\bm b\in \mathbb M_{\bm X}=\widetilde{\mathbb M}_{\bm X} =\{ E[\bm X]: E\in \mathscr{E} \}.$
Note $\vSbep[\langle \bm p,\bm X\rangle]=\vSbep[\langle \bm p,\bm X_i\rangle]=g(\bm p)$ for all $\bm p$ and $i$. It follows that
 $$\bm b\in\mathbb M_{\bm X}=\mathbb M_{\bm X_i}=\widetilde{\mathbb M}_{\bm X_i} =\{ E[\bm X_i]: E\in \mathscr{E} \}, \; i=1,2,\ldots $$
Hence, by Lemma  \ref{moment_v},
$$ \lowCapc\left(\Big|\frac{\bm S_n}{n}-\bm b\Big|\ge \epsilon\right)\le 2\epsilon^{-2}n^{-2}n\Sbep[|\bm X|^2]
\le 2c^2\epsilon^{-2}n^{-1}\to 0. $$
The proof of \eqref{eqlemWLLN2} is completed.

Finally, we show that \eqref{eqlemWLLN3}  is a corollary of \eqref{eqlemWLLN1} and \eqref{eqlemWLLN2}.
Without loss of generality, we assume $\varphi(\bm x)\ge 0$, for otherwise we can replace it by $\varphi+\|\varphi\|$, where $\|\varphi\|=\sup_{\bm x}|\varphi(\bm x)|$.  It   follows from \eqref{eqlemWLLN1}  that
 \begin{align*}
 \limsup_{n\to \infty}\Sbep\left[\varphi\left(\frac{\bm S_n}{n}\right)\right]
 \le & \sup_{\bm x\in \mathbb M_{\bm X}^{\epsilon}}\varphi(\bm x)+ \|\varphi\|\limsup_{n\to \infty}\upCapc\left(\frac{\bm S_n}{n}\not\in\mathbb M_{\bm X}^{\epsilon}\right) \\
 =& \sup_{\bm x\in \mathbb M_{\bm X}^{\epsilon}}\varphi(\bm x)\to \sup_{\bm x\in \mathbb M_{\bm X}}\varphi(\bm x)\; \text{ as } \epsilon\to 0.
 \end{align*}
Now suppose $b\in \mathbb M_{\bm X}$.
By \eqref{eqlemWLLN2},
\begin{align*}
 & \liminf_{n\to \infty}\Sbep\left[\varphi\left(\frac{\bm S_n}{n}\right)\right]\ge
 \liminf_{n\to \infty}\Sbep\left[\varphi\left(\frac{\bm S_n}{n}\right)I\left\{\left|\frac{\bm S_n}{n}-\bm b \right|<\epsilon\right\}\right]\\
 \ge & \inf_{\bm x:|\bm x-\bm b|<\epsilon}\varphi(\bm x)\liminf_{n\to \infty}\upCapc\left(\Big|\frac{\bm S_n}{n}-\bm b\Big|<\epsilon\right) =\inf_{\bm x:|\bm x-\bm b|<\epsilon}\varphi(\bm x)\to  \varphi(\bm b)
\end{align*}
as $\epsilon\to 0$.  By the arbitrariness of $\bm b\in \mathbb M_{\bm X}$,
 $$
 \liminf_{n\to \infty}\Sbep\left[\varphi\left(\frac{\bm S_n}{n}\right)\right]
 \ge  \sup_{\bm b\in \mathbb M_{\bm X}}\varphi(\bm b).
$$
The proof of \eqref{eqlemWLLN3} is completed. $\Box$

\bigskip

{\bf Proof of Theorem \ref{thLLNv}}.  Let $Q$ be a countable subset of $\mathbb R^d$ which is dense in $\mathbb R^d$. Then
 \begin{align*}
  \outCapc\left(C\left\{\frac{\bm S_n}{n}\right\}\not\subset \mathbb M_{\bm X} \right)
 =&\outCapc\left(\bigcup_{\bm p\in\mathbb R^d} \left\{\limsup_{n\to \infty}\frac{\langle\bm p,\bm S_n\rangle}{n}>\vSbep[\langle\bm p,\bm X_1\rangle]\right\} \right) \\
  =&\outCapc\left(\bigcup_{\bm p\in Q} \left\{\limsup_{n\to \infty}\frac{\langle\bm p,\bm S_n\rangle}{n}>\vSbep[\langle\bm p,\bm X_1\rangle]\right\} \right) \\
  \le &
 \sum_{\bm p\in Q} \outCapc\left(   \limsup_{n\to \infty}\frac{\langle\bm p,\bm S_n\rangle}{n}>\vSbep[\langle\bm p,\bm X_1\rangle]  \right)=0
 \end{align*}
by \eqref{eqthLLN1.2}. And so, \eqref{eqthLLNv.2} is proved.

For \eqref{eqthLLNv.3}, is is sufficient to show that
 \begin{equation}\label{eqproofLLNv.10}
  \Capc^{\mathscr{P}} \left(C\left\{\frac{\tilde{\bm S}_n}{n}\right\}\supset \mathbb M_{\bm X} \right)=1.
  \end{equation}
For any $\bm b\in \mathbb M_{\bm X}$ and $\epsilon>0$, by Proposition \ref{lemWLLN}, we have
\begin{equation}\label{eqproofLLNv.11}
   \lim_{n\to\infty}\upCapc\left(\left|\frac{ \bm S_n}{n}-\bm b\right|\le \epsilon\right)=1.
\end{equation}
Let $\Theta=\{\bm b_1,\bm b_2,\ldots\}$ be a countable subset of $\mathbb M_{\bm X}$ which is dense in $\mathbb M_{\bm X}$.
Let  $\epsilon_k=1/2^k$. By \eqref{eqproofLLNv.11},  there exists a sequence $\{n_k\}$ with $n_k\nearrow \infty$, $n_{k-1}/n_k^{1/p}\to 0$ such that
$$ \upCapc\left(\left|\frac{\bm S_{n_k}-\bm S_{n_{k-1}}}{n_k}-\bm b_j\right|\le  \epsilon_k \right)\ge 1/2, \;\; j=1,\ldots k.$$
Denote
$$A_{k,j} =\begin{cases}
\left\{\left|\frac{\bm S_{n_k}-\bm S_{n_{k-1}}}{n_k}-\bm b_j\right|\le \epsilon_k \right\}, &j=1,2,\ldots,k \\
\emptyset, & j>k.
\end{cases}
 $$
 Then  $$ \sum_{k=1}^{\infty}\upCapc\left((A_{k,j}\right)=\sum_{k=j+1}^{\infty}\upCapc\left(A_{k,j}\right)=\infty, \;\; j=1,2,\ldots.$$
 Note that $A_{k,j}$s are close sets of $\bm X=(X_1,X_2\ldots)$.
 By the Borel-Cantelli Lemma (Lemma \ref{lemBC2} (iii)),
 $$\Capc^{\mathscr{P}}\left(\bigcap_{j=1}^{\infty} \big\{A_{k,j}\;\; i.o. \big\}\right) =1. $$
Note on the event $A=\bigcap_{j=1}^{\infty} \big\{ A_{k,j} \;\; i.o. \big\}$ and $B=\left\{C\left\{\frac{ \bm S_n}{n}\right\}\subset \mathbb M_{\bm X}  \right\}$, we have
\begin{align*}
&\liminf_n \left|\frac{\bm S_n}{n}-\bm b_j\right|\le
\liminf_k \left|\frac{\bm S_{n_k}}{n_k}-\bm b_j\right| \\
= &  \liminf_k \left|\frac{\bm S_{n_k}-\bm S_{n_{k-1}}}{n_k}-\bm b_j\right|=0,
\;\; \text{ for all } \bm b_j\in \Theta.
\end{align*}
Note that $\Theta$ is dense in $\mathbb M_{\bm X}$. It follows that on $A$ and $B$,
\begin{align*}
 \liminf_n \left|\frac{\bm S_n}{n}-\bm b \right|=0,
\;\; \text{ for all } \bm b\in \mathbb M_{\bm X}.
\end{align*}
On the other hand, $\Capc^{\mathscr{P}}(B^c)=0$ by \eqref{eqthLLNv.2}. So, $\Capc^{\mathscr{P}}(AB)\ge \Capc^{\mathscr{P}}(A)-\Capc^{\mathscr{P}}(B^c)=1$.
It follows that
$$ \Capc^{\mathscr{P}}\left(\liminf_n \left|\frac{\bm S_n}{n}-\bm b \right|=0
\;  \text{ for all } \bm b\in \mathbb M_{\bm X}\right)=1. $$
Hence, \eqref{eqproofLLNv.10} is proved.

 Finally, we consider \eqref{eqcorLLNv.4}. Let $\bm Y_i=\bm X_i^{(i)}$, $\bm T_n=\sum_{i=1}^n \bm Y_i$, where $\bm X^{(c)}=(X_1^{(c)},\ldots,X_d^{(c)})$ for $\bm X=(X_1,\ldots, X_d)$.  Then
\begin{equation} \label{eqproofLLNv.6} \sum_{n=1}^{\infty}\upCapc(|\bm X_n|>n )\le  \sum_{n=1}^{\infty}\upCapc(|\bm X|>n/2)<\infty,
\end{equation}
\begin{equation} \label{eqproofLLNv.7} \frac{\sum_{i=1}^n \vSbep[|\bm X_i-\bm Y_i|]}{n}\le \frac{\sum_{i=1}^n \sum_{j=1}^d\vSbep[(|X_{i,j}|-i)^+]}{n}\to 0
\end{equation}
and
\begin{equation} \label{eqproofLLNv.8} \sum_{i=1}^{\infty}\frac{\Sbep[|\bm Y_i|^2]}{i^2}\le \sum_{i=1}^{\infty}\frac{\Sbep[|\bm X_1|^2\wedge (di)^2]}{i^2}<\infty,
\end{equation}
 by Lemma \ref{lem3.6}.

  When $\bm b\not\in \mathbb M_{\bm X}$,  \eqref{eqcorLLNv.4} is obvious by  \eqref{eqthLLNv.2}. Suppose
 $$ \bm b\in  \mathbb M_{\bm X}=\widetilde{\mathbb M}_{\bm X_i} =\{ E[\bm X_i]: E\in \mathscr{E} \}, \; i=1,2,\ldots $$
There exists $E_i\in \mathscr{E}$ such that $\bm b=E_i[\bm X_i]$. Note that each $\bm Y_n$ is tight.  For  linear operators $E_i$ and the sequence $\{\bm Y_n;n\ge 1\}$, by Proposition \ref{lem3.4} there exist a copy $\{\tilde{\bm Y}_n;n\ge 1\}$ on $(\widetilde{\Omega},\widetilde{\mathscr{H}},\widetilde{\mathbb E})$ and a probability measure $Q$ on $\sigma(\tilde{\bm Y}_1,\tilde{\bm Y}_2,\ldots)$ such that $\{\tilde{\bm Y}_n;n\ge 1\}$ is a sequence of independent random variables under $Q$,
\begin{equation} \label{eqprooflemWLLN.3}
Q\left[\varphi(\tilde{\bm Y}_i)\right]=E_i\left[\varphi(\bm Y_i)\right]\; \text{ for all } \varphi\in C_{b,Lip}(\mathbb R^d),
\end{equation}
\begin{equation} \label{eqprooflemWLLN.4}
Q\left[\varphi(\tilde{\bm Y}_1,\ldots,\tilde{\bm Y}_p)\right]\le \Sbep\left[\varphi(\bm Y_1,\ldots,\bm Y_p)\right]\; \text{ for all } \varphi\in C_{b,Lip}(\mathbb R^{d\times p})
\end{equation}
and
\begin{equation} \label{eqprooflemWLLN.5}
\widetilde{v}\left( B\right)\le Q\left( B\right)\le \widetilde{V}\left( B\right) \; \text{ for all } B\in\mathscr{B}(\tilde{\bm Y}_1,\tilde{\bm Y}_2,\ldots).
\end{equation}
 Similar to \eqref{eqExpQ} and \eqref{eqVarQ}, we have $Q[\tilde{\bm Y}_i]=E_i[\bm Y_i]$ and $Q[|\tilde{\bm Y}_i|^2]\le \Sbep[|\bm Y_i|^2]$ . Then
\begin{equation} \label{eqproofLLNv.17} \frac{1}{n}\sum_{i=1}^n|Q[\tilde{\bm Y}_i]-\bm b|= \frac{1}{n}\sum_{i=1}^n|E_i[\bm Y_i]-E_i[\bm X_i]|
\le \frac{1}{n}\sum_{i=1}^n\vSbep[|\bm Y_i-\bm X_i|]\to 0
\end{equation}
by \eqref{eqprooflemWLLN.4} and \eqref{eqproofLLNv.7}, and
$$ \sum_{i=1}^{\infty}\frac{Q[|\tilde{\bm Y}_i|^2]}{i^2}
\le \sum_{i=1}^{\infty}\frac{\Sbep[|\bm Y_i|^2]}{i^2} <\infty,
 $$
 by \eqref{eqproofLLNv.8} and \eqref{eqprooflemWLLN.4}. With the same arguments as showing \eqref{eqconvergenceforBC} we have
 $$\sum_{k=1}^{\infty}\cCapc^{\mathscr{P}} \left(  \frac{\max_{n_k+1\le n\le n_{k+1}}|\sum_{i=n_k+1}^n(\bm Y_i-Q[\tilde{\bm Y}_i])|}{n_k}>2\epsilon_k\right)<\infty, $$
where $n_k=2^k$, $\epsilon_k\searrow 0$,  which, similar to  \eqref{eqproofcor1}, implies
 $$\Capc^{\mathscr{P}} \left(\lim_{n\to\infty}\frac{\sum_{i=1}^n(\bm Y_i-Q[\tilde{\bm Y}_i])}{n}=0\right)=1. $$
  On the other hand, by \eqref{eqproofLLNv.6} and the Borel-Cantelli lemma, we have
 $\outCapc(\bm X_n\ne \bm Y_n \; \; i.o.)=0$. It follows that
 $$  \Capc^{\mathscr{P}}\left(\lim_{n\to\infty}\frac{\bm S_n}{n}=\bm b\right)=1. $$
\eqref{eqcorLLNv.4} is proved.
$\Box$

%\Acknowledgements{This work was Supported by grants from the NSF of China (Grant No.11731012,12031005),   Ten Thousands Talents Plan of Zhejiang Province (Grant No. 2018R52042) and the Fundamental Research Funds for the Central Universities.}


\end{document}